\documentclass[a4paper, 12pt]{article}
\usepackage{latexsym}
\usepackage[T1]{fontenc}
\usepackage{amsmath}
\usepackage{amssymb}
\usepackage{graphicx,graphpap}
\usepackage{color}
\usepackage{dlfltxbcodetips}
\usepackage{dsfont}
\usepackage{kbordermatrix}
\newtheorem{defn}{Definition}[section]
\newtheorem{lem}[defn]{Lemma}
\newtheorem{prop}[defn]{Proposition}

\newtheorem{cor}[defn]{Corollary}
\newtheorem{rem}[defn]{Remark}
\newcommand{\vc}[1]{\boldsymbol{#1}}
\newcommand{\p}{\textrm{P}}

\newcommand{\qed}{\hfill $\square$}
\newcommand{\ma}[1]{{\bf{#1}}}

\usepackage{pifont}

\begin{document}

\title{How old is this bird?\\ The age distribution under some phase sampling schemes}
\author{Sophie Hautphenne\footnote{University of Melbourne}, Melanie Massaro\footnote{Charles Sturt University} and Peter Taylor\footnote{University of Melbourne}}

\maketitle

\begin{abstract} In this paper, we use a finite-state continuous-time Markov chain with one absorbing state to model an individual's lifetime. Under this model, the time of death follows a phase-type distribution, and the transient states of the Markov chain are known as \textit{phases}. We then attempt to provide an answer to the simple question ``\textit{What is the conditional age distribution of the individual, given its current phase}''? We show that the answer depends on how we interpret the question, and in particular, on the phase observation scheme under consideration. We then apply our results to the computation of the age pyramid for the endangered Chatham Island black robin \textit{Petroica traversi} during years of intensive conservation efforts in 1980-1989.
\end{abstract}

\section{Introduction}
\label{intro}
A random variable has a phase-type (PH) distribution if it corresponds to the time until absorption of a transient Markov chain with one absorbing state. PH distributions, introduced in the early 1980's by Neuts \cite[Chapter 2]{neuts1981matrix}, form a class of distributions with considerable modelling versatility, which results from attractive probabilistic properties. The set of PH distributions is closed under convolutions and finite mixtures, and is dense in the class of all distributions with non-negative support. PH distributions have therefore been extensively used in practice, in particular for modelling lifetime distributions, see for instance \cite{aalen1995phase},  \cite{gavrilov1991biology}, and \cite{lin2007markov}.

The question addressed in this paper arose initially when modelling the lifetime reproductive success of the black robin \textit{Petroica traversi}, which
 is an endangered songbird species endemic to the Chatham Islands, an isolated archipelago located 800km east of New Zealand. By 1980, the population of black robins had declined to five birds, including a single successful breeding pair, on Mangere Island ~\cite{butlerblack}. Through intensive conservation efforts in 1980-1989 by the New Zealand Wildlife Service (now the Department of Conservation), the population recovered to 93 birds by spring 1990 \cite{kennedy2014severe}. Over the next decade (1990-1998), the population was closely monitored, but without human intervention. Nevertheless the population continued to grow rapidly to 197 adults by 1998, but after this period, the population growth slowed considerably and it only reached 239 adults in 2011 \cite{massaro2013nest}. 


 The black robin population is modeled in a parallel study\footnote{\label{ff}S. Hautphenne, M. Massaro, E. S. Kennedy, and R. Sainudiin. Modelling of the Chatham Island black robin \textit{Petroica traversi} populations using branching processes: Informed management strategies for reintroduction of endangered species. In preparation} using a special class of branching processes called the Markovian binary tree (MBT), in which an underlying transient Markov chain controls the reproduction and death events of each individual in the population. A direct consequence of this model is that each bird lives for a random time which has a PH distribution, in which it progresses through states (also called phases) of a continuous-time Markov chain and dies when the chain moves to an absorbing state. In this application, the phases do not have any particular physical interpretation, their role is to increase the accuracy and realism of the MBT model, as opposed to the simplest linear birth-and-death model. In other real-world applications, the phases may have a physical meaning, such as in \cite{lin2007markov} where they model physiological ages, which can be interpreted as relative health indices, as opposed to chronological age.

MBTs have proved to be powerful stochastic models in population biology and demography \cite{hautphenne2012markovian}. Having fitted an MBT to real data, we can calculate properties of the population, such as the probability that it will become extinct in some time interval $[0, t ]$, and the distribution of the population size at time $t$. In particular, the model allows us to compute the expected number of birds in phase $j$ at time $t$, and the asymptotic frequency of phase $j$ in the population. However, the latter two quantities may not have significance for biologists, who are likely to be interested in age-specific, rather than phase-specific, properties of individuals. We therefore need to be able to translate information about the phase distribution into information about the age distribution; that is, we need to answer the question ``\textit{what is the age distribution of a bird, given its phase}''? 

Note that the reverse exercise of translating information about age into information about phase is much easier since the distribution of the phase at any given age is well known. The main difference is that age is deterministic, while phase is random and an individual stays in a given phase for an exponentially distributed amount of time. 

In fact, defining the event that an individual is in phase $j$ is already not trivial. It is necessary to describe in more detail how the individual is sampled. We suggest three sampling schemes according to which an observer looks at the phase of an individual: \begin{itemize}\item first, we assume that individuals are born according to a Poisson process which started infinitely far in the past, and one observes the phase of a randomly selected individual at time 0;
\item second, no assumption is made on the birth process, and the observations of individuals occur according to a Poisson process, in which case we allow a single or multiple observation(s);
\item third, a single observation occurs at a uniformly distributed random time within some time window.
\end{itemize}
 For each observation scheme, we compute the conditional age distribution given the observed phase, as well as related quantities.  We show that the age distribution conditional on a single rare Poisson phase observation coincides with the age distribution conditional on a single rare uniform phase observation. Moreover, this age-distribution also corresponds to the age distribution of a randomly selected bird in a given phase at time 0 in the process where individuals are born according to a Poisson process.

The questions addressed in this paper, and their proposed answers, are not restricted to the context of PH distributions and ageing processes; they have a wider interest in Markov chain theory. Indeed, if $Q$ is the generator of a continuous-time Markov chain $\{X(t):t\geq 0\}$, then it is well known that $P[X(t)=j|X(0)=i]=[\exp(Qt)]_{ij}$, but the conditional distribution of the time $t$ elapsed since the start of the Markov chain, given that the chain is observed in phase $j$ is much less explored. As we observed above, the nature of this observation event needs to be described carefully.

The paper is organised as follows. In the next section, we provide some background on PH distributions used to model the ageing process of individuals. In Section 3, we compute the conditional age distribution at time 0 in the Poisson birth process, given the observed phase. In section 4, we consider the Poisson phase observation scheme and provide the conditional age distribution with a single or multiple observation(s). In Section 5, we consider a uniform observation scheme, and in Section 6 we discuss the rare observation limit of the results obtained in Sections 4 and 5. Finally, in Section 7, we illustrate our results on a toy example first, and then on the computation of the age pyramid for the black robin population.

\section{The phase-type lifetime distribution}

We assume that the lifetime of an individual is a random variable $L$ which follows a phase-type PH$(\vc\alpha,{\bf{Q}})$ distribution with $m$ transient phases $\{1,2,\ldots,m\}$ and the absorbing phase 0. This PH distribution is parameterised by an $1\times m$ vector $\vc\alpha$ which gives us the initial distribution of the underlying Markov chain, and an $m\times m$ matrix ${\bf Q}$ containing the transition rates between the transient phases. So the assumption is that the lifetime of an individual progresses through phases (which may or may not correspond to some physically-observable characteristics) according to a realisation of the Markov chain, and the individual dies when the chain moves to the absorbing phase 0.

 The PH$(\vc\alpha,{\bf Q})$ distributed random variable $L$ has a density and a distribution function respectively given by
\begin{eqnarray}\nonumber f_L(x)&= &\vc \alpha e^{{\bf Q}x}\vc q_0\\\label{disPH}F_L(x)=\p[L\leq  x]&=&1- \vc \alpha e^{{\bf Q}x}\vc 1,\end{eqnarray}where $\vc q_0=(-{\bf Q})\vc 1$ is the absorption rate vector, and $\vc1$ is a column vector of ones. Let $\varphi(x)$ denote the phase of the individual at age $x$, and let $j$ be any transient phase.
Another basic result on PH random variables tells us that the probability that an individual is in phase $j$ when its age
is $x$ is
$$\p[\mbox{phase}=j\,|\,\mbox{age}=x]=\p[\varphi(x)=j]=(\vc\alpha e^{{\bf Q}x})_j=\vc\alpha e^{{\bf Q}x}\vc e_j,$$ where $\vc e_j$ is the $j$th unit (column) vector. 
Our question is the reverse of this:
``\textit{If we observe an individual in phase j, what can we say about its age?}''
Bayes' Theorem gives us
$$\p[\mbox{age}\leq x\,|\,\mbox{phase}=j] = \dfrac{\p[\mbox{phase}=j\,|\,\mbox{age}\leq x]\,\p[\mbox{age}\leq x]}{\p[\mbox{phase}=j]}.$$ The problem is that we do not yet have anything in the model to make sense of $\p[\mbox{age}\leq x]$ and $\p[\mbox{phase}=j]$. Our aim is  in this paper is to find satisfying answers to the above question.

\section{A Poisson birth process}

One assumption that we might make is that individuals have been born at the epochs of a Poisson process with parameter $\beta$ over the time interval $(-\infty, 0)$, and that we observe the phase of a single individual randomly taken from the population at time 0. Let $F_j(s)$ be the probability that the age of the randomly-selected individual in phase $j$ at time 0 is smaller than $s$.

\begin{lem} \label{pbirth}For any $s\geq 0$,
\begin{equation}\label{Fjs}F_j(s)=1-\dfrac{\vc\alpha e^{{\bf Q}s}(-{\bf Q})^{-1}\vc e_j}{\vc\alpha (-{\bf Q})^{-1}\vc e_j}.\end{equation}
\end{lem}
\noindent \textit{Proof.}
For any $T > 0$,
\begin{itemize}
\item the number $N_T$ of individuals born in the time interval
$[-T, 0)$ has a Poisson distribution with parameter $\beta T$,
\item  conditional on $N_T = n$, the birthtimes of the $n$ individuals are uniformly and independently distributed on the interval $[-T, 0)$,
\item an individual born at time $u \in [-T, 0)$ will be alive at time 0 with probability $\vc \alpha e^{{\bf Q}(-u)}\vc1$ by \eqref{disPH}.
\end{itemize}
So, an individual born in the interval $[-T, 0)$ will be alive at time 0 with probability
$(1/T)\int_{-T}^0 \vc \alpha e^{{\bf Q}(-u)}\vc1\,du.$ It will be alive and in phase $j$ at time 0 with probability
$(1/T)\int_{-T}^0\vc \alpha e^{{\bf Q}(-u)}\vc e_j\,du,$ and it will be alive, in phase $j$
and older than $s$ at time 0 with probability
$(1/T)\int_{-T}^s \vc \alpha e^{{\bf Q}(-u)}\vc e_j\,du.$ So the probability that an individual in phase $j$ at time 0 is older
than $s$ is 
$$\dfrac{\int_{-T}^s \vc \alpha e^{{\bf Q}(-u)}\vc e_j\,du}{\int_{-T}^0 \vc \alpha e^{{\bf Q}(-u)}\vc e_j\,du}.$$ The fact that any PH distribution has a finite mean allows us to let $T\rightarrow\infty $ and, changing the variable of integration, we arrive at the conclusion that the probability that a randomly-selected individual in phase $j$ at time 0 is older than $s$ is
$$\dfrac{\int_{s}^\infty\vc \alpha e^{{\bf Q}u}\vc e_j\,du}{\int_{0}^\infty \vc \alpha e^{{\bf Q}u}\vc e_j\,du}=\dfrac{\vc\alpha e^{{\bf Q}s}(-{\bf Q})^{-1}\vc e_j}{\vc\alpha (-{\bf Q})^{-1}\vc e_j},$$which completes the proof.\qed
\medskip

Note that this model corresponds to an $M/PH/\infty$ queue, and \eqref{Fjs} gives the distribution of the age of a randomly selected individual in steady state. Also observe that $F_j(s)$ does not depend on the rate $\beta$ of the Poisson birth process.

The above analysis is nice. However the birth process in an MBT is not Poisson and, more generally, there is no reason to believe that a Poisson process is a good model for births.
In the next three sections, we shall follow an alternative approach: without making any assumption on the birth process, we look at just a single individual and proceed by explicitly putting the observation process into the model. 

For further use, we denote the age of the individual at the time of observation as $A_o$ and the observed phase as $\varphi_o$. We are therefore interested in computing \begin{equation} \label{Ao}\p[A_o\leq s\,|\,\varphi_o=j].\end{equation}


\section{The Poisson observation scheme}

In this section, we assume that, following its birth, the phase of an individual is observed according to a Poisson process with rate $\gamma$. We first compute the age distribution given the phase at the first observation time. Then we generalise our results to the age distribution at the time of the last observation, given that the observer records the phases at $k\geq 2$ successive time events of the Poisson process.

\subsection{Single observation}
By the properties of Poisson processes, the rate at which an individual is observed when it is in phase $j$ is $\gamma$, for any $j$. A slight modification of the underlying phase process then allows us to compute the conditional age distribution of the individual at the first observation time, given that the individual is in phase $1\leq j\leq m$ at that time. 
It suffices to add $m$ absorbing phases $1',2',\ldots,m'$ (one per transient phase), to the process so that phase $j'$ is reached when the individual is observed in phase $j$. The initial distribution $\vc\alpha$ stays unchanged, but now the transition rate matrix becomes
$$\bf T(\gamma)={\bf Q}-\gamma{\bf  I},$$
and there are $m + 1$ absorption rate vectors $$\vc q_0, \vc{t}_{1'}(\gamma) ,\ldots,\vc{t}_{m'}(\gamma),$$ where $\vc q_0=-{\bf Q}\vc1$ records the rates of absorption into phase 0 (corresponding to the death of the individual), and for $1'\leq j'\leq m'$, $\vc{t}_{j'}(\gamma)=\gamma\vc e_j$ records the rates of absorption into phase $j'$ (corresponding to the observation of the individual in phase $j$). For the sake of clarity of the presentation, we shall drop the dependence on $\gamma$ in ${\bf T}(\gamma)$ and $\vc{t}_{j'}(\gamma)$ and use the simpler notation ${\bf T}$ and $\vc{t}_{j'}$ in the sequel.


For any phase $j$ (transient or absorbing), let $$B(j)=\inf\{t\geq 0: \varphi(t)=j\}$$ be first time the individual enters phase $j$, with $B(j) = \infty$ if the individual never enters phase $j$. Then, with probability one, precisely one of the random variables
$$\{B(0), B(1'), B(2'), . . . , B(m')\}$$ is finite, and the age distribution conditional on the observed
phase being $j$ can be rewritten as \begin{equation}\label{q0}\p[A_o\leq s\,|\,\varphi_o=j]=\p[B(j')\leq s \,|\,B(j')<\infty].\end{equation} 
Based on this observation, the next proposition provides an expression for the conditional age distribution.


 \begin{prop}\label{poissingle}The age distribution of the individual at the first observation (event) time of a Poisson$(\gamma)$ process, conditional on the observed phase being $j$, is given by \begin{equation}\label{q0sol}\p[A_o\leq s\,|\,\varphi_o=j]=1-\dfrac{ \vc\alpha \exp({\bf T} s)\,(-{\bf T})^{-1} \vc{e}_{j}}{\vc\alpha (-{\bf T})^{-1} \vc{e}_j}.\end{equation}
 \end{prop}
 
\noindent \textit{Proof.} We have
 \begin{eqnarray*}\p[B(j')\leq s \,|\,B(j')<\infty]&=&1-\p[B(j')> s \,|\,B(j')<\infty]\\&=&1-\dfrac{\p[s<B(j') <\infty]}{\p[B(j') <\infty]}.
\end{eqnarray*}On the one hand,
\begin{eqnarray*}\p[s<B(j') <\infty]&=&\int_s^\infty
\vc\alpha \exp({\bf T}u)\,\vc{t}_{j'} du\\&=&\vc\alpha \exp({\bf T} s)\,(-{\bf T})^{-1} \gamma \vc{e}_{j}, \end{eqnarray*} and on the other hand, since $\p[B(j')>0]=1$, \begin{equation}\label{ppo}\p[B(j') <\infty]=\p[0<B(j') <\infty]=\vc \alpha \,(-{\bf T})^{-1} \gamma\vc{e}_{j},\end{equation} which, with \eqref{q0}, completes the proof.
\qed
\medskip

\begin{figure}[t] 
	\begin{center}
	\begin{picture}(300,150)(0,10)
	\includegraphics[angle=0,width=10cm]{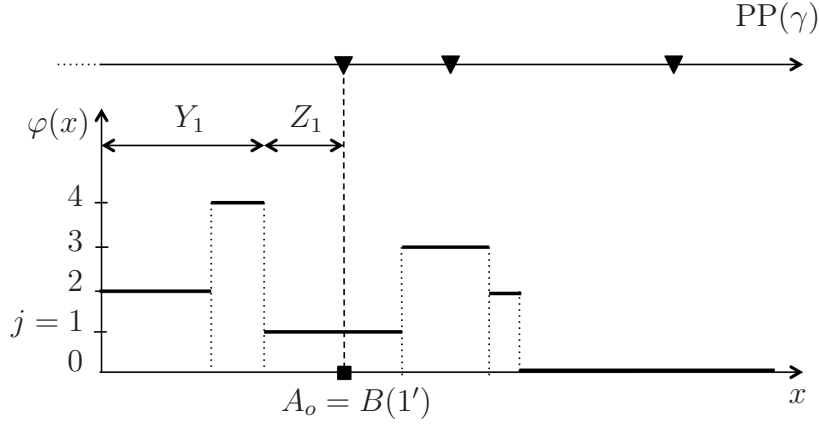} 
	\put(-280,5){0}
	\put(-301,20){$j=1$}\put(-280,35){2}\put(-280,50){3}\put(-280,67){4}\put(-240,97){$Y_1$}\put(-197,97){$Z_1$}\put(-200,-10){$A_o=B(1')$}\put(-10,-7){$x$}\put(-295,95){$\varphi(x)$}\put(-30,135){PP$(\gamma)$}
	\end{picture}\vspace*{1cm}
	\caption{\label{figpath1} A Possible trajectory of the phase process until absorption. The observation process is a Poisson process with rate $\gamma$ (denoted as PP$(\gamma)$), where the $ \blacktriangledown$ symbols represent the events, and the phase is $j=1$ at the first observation event.}
	\end{center}\end{figure}
	

The age at the observation time, conditional on the observed phase being $j$, can be written as the sum of two random variables, $$A_o=Y_j+Z_j,$$ where $Y_j$ denotes the last time that the Markov chain entered phase $j$ before observation, and $Z_j$ denotes the sojourn time in phase $j$ between this time and observation, both random variables being conditionally independent given $\varphi_o=j$. The random variables $A_o, Y_j, $ and $Z_j$ are illustrated in Figure \ref{figpath1}.  Besides purely theoretical interest, the distribution of $Y_j$ and $Z_j$ may have practical interest when the phases have some physical interpretation (such as physiological ages), and an observer who sees an individual in phase $j$ is interested in knowing the chronological age at which the individual entered that particular phase ($Y_j$), or for how long he/she has been in that phase ($Z_j$).
The respective distributions of $Y_j$ and $Z_j$ are computed in the next two propositions.

\begin{prop} \label{pyjp}The conditional distribution of $Y_j$, given $\varphi_o=j$, has a point mass at zero given by
\begin{eqnarray}\label{Yj0a} \p[Y_j=0\,|\,\varphi_o=j]&=&\dfrac{\alpha_j}{(\lambda_j+\gamma)\,\vc\alpha (-{\bf{T}})^{-1}\vc{e}_j},\end{eqnarray} and for $y>0$,
\begin{eqnarray} \label{Yja}\p[Y_j\leq y\,|\,\varphi_o=j]&=&1-\dfrac{\vc\alpha e^{{\bf T}y}(-{\bf T})^{-1}({\bf Q}+\lambda_j \,{\bf I})\vc e_j}{(\lambda_j+\gamma)\,\vc\alpha (-{\bf T})^{-1}\vc e_j},\end{eqnarray}
where $\lambda_j=-Q_{jj}$.
\end{prop}
\noindent \textit{Proof.} First, recall from \eqref{ppo} that $$\p[\varphi_o=j]=\p[B(j') <\infty]=\vc\alpha \,(-{\bf T})^{-1} \gamma\vc{e}_{j}.$$ Let $T_o$ be the time at which the individual is observed, the clock being set at the individual's birth time. Thanks to the memoryless property of exponential random variables, $T_o$ has the same distribution as the interarrival time in the Poisson observation process, that is, $T_o$ is exponentially distributed with parameter $\gamma$. Then, for any $y>0$, by conditioning on the value of $T_o$, we have
\begin{eqnarray*}\lefteqn{\p[Y_j\in[y,y+dy], \,\varphi_o=j]}\\&=&\int_0^\infty\p[Y_j\in[y,y+dy], \,\varphi(u)=j\,|\,T_o\in[u,u+du]]\,\gamma e^{-\gamma u}du\\&=& \int_y^\infty\sum_{k\neq j}(\vc\alpha e^{\ma Q y})_k Q_{kj} e^{-\lambda_j (u-y)}\,dy\,\gamma e^{-\gamma u}du.
\end{eqnarray*} Next, observe that since $Q_{jj}=-\lambda_j,$ we have $$\sum_{k\neq j}(\vc\alpha e^{\ma Q y})_k Q_{kj} =\vc\alpha e^{\ma Q y} (\ma Q+\lambda_j\,\ma I)\vc e_j.$$ As 
$$ \int_y^\infty e^{-\lambda_j (u-y)}\,\gamma e^{-\gamma u}du=\dfrac{\gamma e^{-\gamma y }}{\lambda_j+\gamma}, $$and $\ma T=\ma Q-\gamma\ma I$, we have
\begin{equation}\label{aaa}\p[Y_j\in[y,y+dy], \,\varphi_o=j]=\dfrac{\vc\alpha e^{\ma T y} (\ma Q+\lambda_j\,\ma I)\vc e_j\gamma\,dy}{\lambda_j+\gamma}.\end{equation}
Similarly, $$ \p[Y_j=0, \,\varphi_o=j]=\int_0^\infty \alpha_j e^{-\lambda_j u}\,\gamma e^{-\gamma u}du=\dfrac{\alpha_j \gamma}{\lambda_j+\gamma}, $$ which leads to \eqref{Yj0a}. Finally, \eqref{Yja} follows from \eqref{Yj0a} and \eqref{aaa} since$$\p[Y_j\leq y\,|\,\varphi_o=j]=\p[Y_j=0\,|\,\varphi_o=j]+\int_0 ^y \p[Y_j\in[u,u+du], \,\varphi_o=j],$$ for any $y>0$. 
\qed
\medskip

\begin{prop}\label{pzjp} The conditional distribution of $Z_j$, given $\varphi_o=j$, is exponential with parameter $\gamma+\lambda_j$.
\end{prop}
\noindent \textit{Proof.} Let $S_j$ denote the sojourn time of the underlying Markov chain in phase $j$, and let $E_{j\rightarrow j'}$ denote the event that upon leaving phase $j$, the chain moves to phase $j'$. We have $\p[E_{j\rightarrow j'}]=\gamma/(\lambda_j+\gamma),$ and $$\p[S_j > z, E_{j\rightarrow j'}] =  e^{-(\lambda_j + \gamma) z} \frac{ \gamma} {\lambda_j + \gamma},$$therefore$$\p[Z_j>z\,|\, \varphi_o=j]=\dfrac{\p[S_j > z, E_{j\rightarrow j'}]}{\p[E_{j\rightarrow j'}]}= e^{-(\lambda_j + \gamma) z},$$which proves the statement of the proposition.
\qed

\subsection{Multiple observations} We now assume that the observer makes $k\geq 2$ observations of an individual according to a Poisson process with rate $\gamma$. We further assume that the individual is still living at the time of the last ($k$th) observation. We shall compute the age distribution at the time of the $k$th observation, given the sequence of phases observed at the observation times. 

For that purpose, we consider the same process as in the single observation case, with absorbing phases $0,1',\ldots,m'$. Every observation event corresponds to a phase absorption in one of the phases $1',\ldots,m'$. After a phase absorption in $j'$ (which corresponds to an observation of phase $j$), the process instantaneously starts again in phase $j$, that is, with the initial distribution vector $\vc e_j^\top$, until the next absorption event. 

In order to properly define the quantities of interest, we need to redefine the random variables $B(j)$ as follows: for any initial phase distribution $\vc\theta$ and for any phase $j$, 
\begin{equation}\label{agemult}B_{\vc\theta}(j)=\inf\{t\geq 0: \varphi(t)=j, \varphi(0)\sim \vc\theta\}\end{equation}is the first time the process reaches phase $j$, given that the initial phase follows the distribution $\vc\theta$. For the sake of clarity, we shall write $B_\ell(j)$ instead of $B_{\vc e_\ell^\top}(j)$ when the process starts in phase $\ell$ with probability one. Let $j_1,\ldots,j_k$ be the $k$ successive observed phases. The age of the individual at the last observation time conditional on the observed phases, denoted by $A_o(j_1,\ldots,j_k)$ (or by $A_o$ when there is no confusion), is then given by $$A_o(j_1,\ldots,j_k)=B_{\vc\alpha}(j_1')+\sum_{i=2}^k B_{{j_{(i-1)}}}(j_i').$$ The age distribution at the time of the last observation conditional on the sequence of observed phases can then be written as
\begin{eqnarray*}\p[A_o\leq s\,|\,B_{\vc\alpha}(j_1')<\infty, B_{{j_1}}(j_2')<\infty,\ldots,B_{{j_{k-1}}}(j_k')<\infty].\end{eqnarray*} 
In order to compute this distribution, we need the following lemma, which is a particular case of Theorem~1 in \cite{carbonell2008computing}:
\begin{lem}\label{lemcarb}For $2\leq i\leq k$, define the $m\times m$ matrix \begin{eqnarray*}\lefteqn{{\bf{\mathcal{B}}}_{1i}(s)=}\\&&\hspace{-3mm}\int_{u_1=0}^s\int_{u_2=0}^{s-u_1}\ldots\int_{u_{i-1}=0}^{s-u_1-\ldots-u_{i-2}} e^{\ma A_{11}u_1}\ma A_{12}\,e^{\ma A_{22}u_2}\ma A_{23}\ldots \ma A_{(i-1)i}\, e^{\ma A_{ii}(s-u_1-\ldots-u_{i-1})}d\vc u,\end{eqnarray*}where $s>0$ and $\ma A_{ij}$ are constant $m\times m$ matrices. If the $km\times km$ block-structured matrix $\ma A^{(k)}$ is defined by
$$\ma A^{(k)}=\left[\begin{array}{cccccc}\ma A_{11} & \ma A_{12}&\ma 0&\ma 0 &\ldots & \ma 0\\ \ma 0 &\ma A_{22}&\ma A_{23}&\ma 0&\ldots&\ma 0\\\vdots&&\ddots& & &\vdots\\\ma 0&\ldots & \ma 0&\ma 0&\ma A_{(k-1)(k-1)}&\ma A_{(k-1)k}\\\ma 0&\ldots & \ma 0&\ma 0&\ma 0&\ma A_{kk}\end{array}\right],$$ then\begin{equation}\label{B1i} \mathcal{B}_{1i}(s)=(\vc f_{k,1}^\top\otimes\ma  I_m) e^{\ma A^{(k)}s}(\vc f_{k,i}\otimes \ma I_m),\end{equation}where 
$\vc f_{k,i}$ is a $k\times 1$ unit vector such that $(\vc f_{k,i})_j=\delta_{ij}.$ 
\qed
\end{lem}

For the purpose of computing the conditional distribution of $A_o$, we define the $km\times km$ matrix $\ma A^{(k)}$ for any $k\geq 2$ as \begin{equation}\label{Ak}\ma A^{(k)}=\left[\begin{array}{cccccc} \ma T & \vc e_{j_1'}\vc e_{j_1}^\top&\ma 0&\ma 0 &\ldots &\ma 0\\ \ma 0 & \ma T&\vc e_{j_2'}\vc e_{j_2}^\top&\ma 0&\ldots&\ma 0\\\vdots&&\ddots& & &\vdots\\\ma 0&\ldots & \ma 0&\ma 0&\ma T&\vc e_{j_{k-1}'}\vc e_{j_{k-1}}^\top\\\ma 0&\ldots & \ma 0&\ma 0&\ma 0&\ma T\end{array}\right].\end{equation}

\begin{prop}For an arbitrary $k\geq 2$,  the age distribution of the individual at the $k$th observation time, conditional on the successive observed phases being $j_1,j_2,\ldots,j_k$, is given by
\begin{eqnarray}\label{rap}\p[A_o\leq s\,|\,B_{\vc\alpha}(j_1')<\infty, B_{{j_1}}(j_2')<\infty,\ldots,B_{{j_{k-1}}}(j_k')<\infty]=\dfrac{N_k(s)}{D_k},\end{eqnarray} where 
\begin{eqnarray} \nonumber N_k(s)&=& \vc \alpha(I-e^{\ma Ts})(-\ma T)^{-1}\vc e_{j_1'}\prod_{i=1}^{k-1} \vc e_{j_i}^\top(-\ma T)^{-1}\vc e_{j_{i+1}'}\\\label{Nk}&-&\vc\alpha \sum_{i=2}^k\mathcal{B}_{1i}(s)(- \ma T)^{-1}\vc e_{j_i'}\prod_{\ell=i}^{k-1} \vc e_{j_\ell}^\top(-\ma T)^{-1}\vc e_{j_{\ell+1}'},\\\label{Dk}
D_k&=&\vc \alpha(-\ma T)^{-1}\vc e_{j_1'}\prod_{i=1}^{k-1} \vc e_{j_i}^\top(-\ma T)^{-1}\vc e_{j_{i+1}'},\end{eqnarray} where $\mathcal{B}_{1i}(s)$ is defined in \eqref{B1i} and $\ma A^{(k)}$ is given in \eqref{Ak}.
\end{prop}

\noindent \textit{Proof.} We have
\begin{eqnarray*}\nonumber\lefteqn{\p[A_o\leq s\,|\,B_{\vc\alpha}(j_1')<\infty, B_{{j_1}}(j_2')<\infty,\ldots,B_{{j_{k-1}}}(j_k')<\infty]}\\&=&
\dfrac{\p[A_o\leq s, \,B_{\vc\alpha}(j_1')<\infty, \ldots,B_{{j_{k-1}}}(j_k')<\infty]}{\p[B_{\vc\alpha}(j_1')<\infty, B_{{j_1}}(j_2')<\infty,\ldots,B_{{j_{k-1}}}(j_k')<\infty]}\\&=:&\dfrac{\bar{N}_k(s)}{\bar{D}_k}.\qquad\end{eqnarray*} 
We shall prove using induction on $k$ that \begin{equation}\label{nk}\dfrac{\bar{N}_k(s)}{\bar{D}_k}= \dfrac{{N}_k(s)}{{D}_k},\end{equation} 
where ${N}_k(s)$ and ${D}_k$ satisfy \eqref{Nk} and \eqref{Dk}, respectively. 
Recall that $\vc t_{j'}=\gamma \vc e_j$ for any absorbing phase $1'\leq j'\leq m'$.
When $k=2$, 
\begin{eqnarray}\nonumber\bar{D}_2&=&\p[B_{\vc\alpha}(j_1')<\infty, B_{{j_1}}(j_2')<\infty]\\\label{D2}&=&\vc\alpha (-\ma T)^{-1} \vc t_{j_1'}\,\vc e_{j_1}^\top (-\ma T)^{-1} \vc t_{j_2'}\\\nonumber&=&\gamma^2 D_2.
\end{eqnarray} Further, by \eqref{agemult} and by conditioning on the value of the absorption times $B_{\vc\alpha}(j_1')$ and $B_{{j_{1}}}(j_2')$, we have 
\begin{eqnarray}\nonumber\bar{N}_2(s)&=& \p[(B_{\vc\alpha}(j_1')+B_{{j_{1}}}(j_2'))\leq s, \,B_{\vc\alpha}(j_1')<\infty, B_{{j_1}}(j_2')<\infty]\\\nonumber&=&\int_{u=0}^s \vc\alpha e^{\ma Tu}\vc t_{j_1'}\int_{v=0}^{s-u} \vc e_{j_1}^\top e^{\ma Tv}\vc t_{j_2'}\, dv \,du\\\nonumber&=&\gamma^2\int_{u=0}^s \vc\alpha e^{\ma Tu}\vc e_{j_1'}\vc e_{j_1}^\top (\ma I-e^{\ma T(s-u)})(-\ma T)^{-1}\vc e_{j_2'} \,du\\\label{N2}&=&\gamma^2\{\vc \alpha (\ma I-e^{\ma Ts})(-\ma T)^{-1}\vc e_{j_1'}\vc e_{j_1}^\top(-\ma T)^{-1}\vc e_{j_2'}-\vc \alpha \mathcal{B}_{12}(s)(-\ma T)^{-1}\vc e_{j_2'}\}\qquad\\\nonumber &=&\gamma^2\,N_2(s),\qquad
\end{eqnarray}where \begin{eqnarray*}\mathcal{B}_{12}(s)&=&\int_{u=0}^s e^{\ma Tu} \vc e_{j_1'}\vc e_{j_1}^\top e^{\ma T(s-u)}du.\end{eqnarray*}Using Lemma \ref{lemcarb}, this matrix integral can be evaluated explicitly by defining the $2m\times 2m$ block-structured matrix $$\ma A^{(2)}=\left[\begin{array}{cc} \ma T & \vc e_{j_1'}\vc e_{j_1}^\top\\ \ma 0 & \ma T\end{array}\right],$$so that 
$$\mathcal{B}_{12}(s)=(\vc f_{2,1}^\top\otimes \ma I_m) e^{\ma A^{(2)}s}(\vc f_{2,2}\otimes \ma I_m).$$ Therefore \eqref{nk} holds for $k=2$.


 We now assume that \eqref{nk} holds for $k$, and we need to prove that is still holds for $k+1$. We can decompose the conditional age at the $(k+1)$st observation, $A_o(j_1,\ldots,j_{k+1})$, into the sum of the random variables $B_{\vc\alpha}(j_{1}')$ and $A_o(j_2,\ldots,j_{k+1})$, which are conditionally independent given $j_1$. Note that $A_o(j_2,\ldots,j_{k+1})$ is now conditional on the phase process starting with initial distribution vector $\vc e_{j_1}^\top$ rather than $\vc\alpha$, and the first observed phase is $j_2$ rather than $j_1$, etc.	
To avoid confusion, we shall use the notation $\hat{A}_o(j_2,\ldots,j_{k+1})$ (or $\hat{A}_o$ for short), $\hat{N}_k(s),\hat{D}_k,\hat{\mathcal{B}}_{1i}(s)$ whenever we will be in that situation.

We use the convolution formula for the sum of the two conditionally independent variables $B_{\vc\alpha}(j_{1}')$ and $\hat{A}_o$, together with the conditional distribution of $B_{\vc\alpha}(j_{1}')$ given in \eqref{q0sol} and the induction assumption, to obtain
\begin{eqnarray*}\nonumber\lefteqn{\p[(B_{\vc\alpha}(j_{1'})+\hat{A}_o)\leq s\,|\,B_{\vc\alpha}(j_1')<\infty, B_{{j_1}}(j_2')<\infty,\ldots,B_{{j_{k}}}(j_{k+1}')<\infty]}\\&=&\int_0^s \dfrac{\vc\alpha e^{ \ma Tu}\vc t_{j_1'}}{\vc\alpha (-\ma T)^{-1}\vc t_{j_1'}}\p[\hat{A}_o\leq s-u|B_{{j_1}}(j_2')<\infty,\ldots,B_{{j_{k}}}(j_{k+1}')<\infty]\,du\\&=& \int_0^s \dfrac{\vc\alpha e^{\ma Tu}\vc e_{j_1'}}{\vc\alpha (-\ma T)^{-1}\vc e_{j_1'}}\dfrac{\hat{N}_k(s-u)}{\hat{D}_k}du.
\end{eqnarray*}We immediately see that the denominator of the above expression, $D_{k+1}:=\vc\alpha (-\ma T)^{-1}\vc e_{j_1'}\hat{D}_k$, corresponds to \eqref{Dk} for $k+1$. It remains to show that the numerator, $N_{k+1}(s):= \int_0^s {\vc\alpha e^{\ma Tu}\vc e_{j_1'}}{\hat{N}_k(s-u)}du$, corresponds to \eqref{Nk} for $k+1$. Using \eqref{Nk} and letting $\vc r_{i,k}=(- \ma T)^{-1}\vc t_{j_{i+1}'}\prod_{\ell=i}^{k-1} \vc e_{j_{\ell+1}}^\top(-\ma T)^{-1}\vc t_{j_{\ell+2}'}$, we have
\begin{eqnarray*}\lefteqn{\int_0^s {\vc\alpha e^{\ma Tu}\vc t_{j_1'}}{\hat{N}_k(s-u)}\,du}\\&=&\int_0^s {\vc\alpha e^{\ma Tu}\vc t_{j_1'}}\,\big\{\vc e_{j_1}^\top(\ma I-e^{\ma T(s-u)})\vc r_{1,k}-\vc e_{j_1}^\top \sum_{i=2}^k\hat{\mathcal{B}}_{1i}(s-u)\vc r_{i,k}\big\}\, du\\&=&\vc\alpha (\ma I-e^{\ma Ts})(-\ma T)^{-1}\vc t_{j_1'} \vc e_{j_1}^\top \vc r_{1,k} -\vc\alpha \mathcal{B}_{12}(s) \vc r_{1,k}\\&&-\vc\alpha \sum_{i=2}^k\int_0^s e^{\ma Tu} \vc t_{j_1'} \vc e_{j_1}^\top \hat{\mathcal{B}}_{1,i}(s-u) \,du \,\vc r_{i,k}.
\end{eqnarray*} Using Lemma \ref{lemcarb} and \eqref{B1i}, we can show that 
$$\int_0^s e^{\ma Tu} \vc t_{j_1'} \vc e_{j_1}^\top \hat{\mathcal{B}}_{1,i}(s-u) \,du = \mathcal{B}_{1,i+1}(s),$$ so that by properly redefining the indices we finally obtain what we need.
\qed
\medskip

 Using \eqref{B1i}, the expressions for $N_k(s)$ and $D_k$ can be rewritten as 
\begin{eqnarray}\label{shortN}N_{k}(s)&=&\vc\alpha \vc u_k-\vc\alpha e^{\ma Ts}\vc u_k-\vc v_k e^{\ma A^{(k)}s}\vc w_k, \\\label{shortD} D_k&=&\vc\alpha\vc u_k\end{eqnarray}where
\begin{eqnarray}\label{uk}\vc u_k&=& (-\ma T)^{-1}\vc e_{j_1'}\prod_{i=1}^{k-1} \vc e_{j_i}^\top(-\ma T)^{-1}\vc e_{j_{i+1}'}\\\label{vk}\vc v_k&=&(\vc f_{k,1}^\top\otimes\vc \alpha)\\\nonumber\vc w_k&=&
 \sum_{i=2}^{k}(\vc f_{k,i}\otimes I_m)(- \ma T)^{-1}\vc e_{j_i'}\prod_{\ell=i}^{k-1} \vc e_{j_\ell}^\top(-\ma T)^{-1}\vc e_{j_{\ell+1}'}.
\end{eqnarray} In the expression for $\vc w_k$, an empty product (when $i=k$) is interpreted as the scalar 1. Note that it is also possible to express $N(k)$ and $D(k)$ recursively as follows: for $k\geq 3$,
\begin{eqnarray*} N_k(s)&=&[N_{k-1}(s)  \vc e_{j_{k-1}}-(\vc f_{k,1}^\top\otimes \vc\alpha) e^{\ma A^{(k)}s}(\vc f_{k,k}\otimes I_m)](-\ma T)^{-1}\vc t_{j_k'} \\
D_k&=&D_{k-1}\vc e_{j_{k-1}}(-\ma T)^{-1}\vc t_{j_k'},\end{eqnarray*}where $N_2(s)$ and $D_2$ are given in \eqref{N2} and \eqref{D2} respectively.
\medskip

We now assume that the individual is still alive at the time of the $k$th observation, but is discovered dead (that is, in phase 0) at the time of the $(k+1)$st observation. We are then interested in the conditional lifetime distribution of the individual, given the sequence of observed phases. Indeed, the lifetime $L$ it is then given by the age at the $k$th observation plus the time until absorption from the last observed phase $j_k$ to phase 0, conditional on this time being less than the time between the $k$th and the $(k+1)$st observation. 
We shall need the following lemma:

\begin{lem}\label{lemph} Let $X\sim$ PH$(\vc\theta, T)$ and $Y\sim Exp(\gamma)$. The conditional distribution of $X$, given that $X\leq Y$, is given by
\begin{equation*}\p[X\leq x\,|\,X\leq Y]=\dfrac{1-\gamma\vc\theta(\gamma\ma I-\ma T)^{-1}\vc 1+\vc \theta e^{(\ma T-\gamma\ma I)x}(\gamma\ma I-\ma T)^{-1}\ma T\vc1}{1-\gamma\,\vc\theta(\gamma\ma I-\ma T)^{-1}\vc 1},
\end{equation*} and the density is given by
\begin{equation*}\label{resid}
f_{X|X\leq Y}(x)=\dfrac{\vc \theta e^{(\ma T-\gamma\ma I)x}(-\ma T)\vc1}{1-\gamma\vc\theta(\gamma\ma I-\ma T)^{-1}\vc 1}.
\end{equation*}
\end{lem}

\noindent \textit{Proof.} The distribution is obtained by conditioning on the value of $Y$. The expression for the conditional density then follows.
\qed
\medskip

The conditional lifetime distribution is computed in the next proposition.
\begin{prop} The lifetime distribution of an individual, conditional on the sequence of observed phases being $j_1,j_2,\ldots, j_{k}, 0$, is given by
\begin{eqnarray}\nonumber\lefteqn{\p[L\leq s\,|\,B_{\vc\alpha}(j_1')<\infty, \ldots,B_{{j_{k-1}}}(j_{k}')<\infty,B_{{j_{k}}}(0)<\infty]}\\\nonumber&&=C_k\left\{ \vc\alpha \vc u_{k}\vc e_{j_{k}}^\top	 [e^{(\ma T-\gamma\ma I)s}-\ma I](\gamma \ma I-\ma T)^{-1}+\vc\alpha\mathcal{I}_{12}(s)+\vc v_k\mathcal{J}_{1k}(s)\right\}\ma T\vc 1, \\\label{lifesol}&&
\end{eqnarray}where $\vc u_{k}$ and  $\vc v_k$ are given in \eqref{uk} and \eqref{vk} respectively, and\begin{eqnarray}\label{Ck}C_k&=&[ \vc\alpha \vc u_{k}(1-\gamma\vc e_{j_{k}}^\top(\gamma \ma I-\ma T)^{-1}\vc 1) ]^{-1},\\\nonumber\mathcal{I}_{12}(s)&=&(\vc f_{2,1}^\top\otimes \ma I_m) e^{\ma B^{(k)}s}(\vc f_{2,2}\otimes \ma I_m),\\ \nonumber\mathcal{J}_{1k}(s)&=&[\ma I_{mk},\ma 0_{mk\times m}]e^{\ma C^{(k)}s}[\ma 0_{m\times mk},\ma I_m]^\top\end{eqnarray} with 
$$\ma B^{(k)}=\left[\begin{array}{cc} \ma T &\vc u_{k}\vc e_{j_{k}}^\top\\ \ma 0 &\ma T-\gamma \ma I\end{array}\right],\quad \ma C^{(k)}=\left[\begin{array}{cc} \ma A^{(k)} &\vc w_{k}\vc e_{j_{k}}^\top\\ \ma 0 &\ma T-\gamma\ma  I\end{array}\right].$$

\end{prop}

\noindent \textit{Proof.}
We have 
\begin{eqnarray}\nonumber\lefteqn{\p[L\leq s\,|\,B_{\vc\alpha}(j_1')<\infty, \ldots,B_{{j_{k-1}}}(j_{k}')<\infty,B_{{j_{k}}}(0)<\infty]}\\\label{life}&=& \p[(B_{\vc\alpha}(j_1')+\sum_{i=2}^{k} B_{{j_{(i-1)}}}(j_i'))+X\leq s\,|\,B_{\vc\alpha}(j_1')<\infty, \ldots,B_{{j_{k}}}(0)<\infty]\phantom{bbbbb}\end{eqnarray} where $X\sim$ PH$(\vc e_{j_{k}}^\top,\ma T)$ and $X$ is taken conditionally on $X\leq Y$, where $Y$ is the interarrival time in the Poisson observation process, which is exponentially distributed with parameter $\gamma$. Since the age at the $k$th observation, given by $B_{\vc\alpha}(j_1')+\sum_{i=2}^{k} B_{{j_{(i-1)}}}(j_i')$, and the residual life time, $X$, are conditionally independent given the phase at the $k$th observation, $j_{k}$, we can use \eqref{rap}--\eqref{Dk}, together with Lemma \ref{lemph} and the convolution formula, to compute the conditional lifetime distribution. In order to simplify the notation we define $C_k$ as in \eqref{Ck}.
 By conditioning on the value of $X$ and using \eqref{shortN}, we then obtain
\begin{eqnarray}\nonumber\lefteqn{\p[L\leq s\,|\,B_{\vc\alpha}(j_1')<\infty, \ldots,B_{{j_{k-1}}}(j_{k}')<\infty,B_{{j_{k}}}(0)<\infty]}\\\nonumber&=&\int_0^s \dfrac{N_k(s-u)}{D_k}\dfrac{\vc e_{j_{k}}^\top e^{(\ma T-\gamma \ma I)u}\, (-\ma T)\vc 1}{(1-\gamma\vc e_{j_{k}}^\top(\gamma I-\ma T)^{-1}\vc 1)} du\\\nonumber&=&C_k\left\{ \vc\alpha \vc u_{k}\vc e_{j_{k}}^\top [\ma I-e^{(\ma T-\gamma \ma I)s}](\gamma \ma I-\ma T)^{-1}-\vc\alpha \int_0^s e^{\ma T(s-u)}\vc u_{k}\vc e_{j_{k}}^\top e^{(\ma T-\gamma \ma I)u} du \right.\\\nonumber&&\left.-\vc v_k\int_0^s e^{\ma A^{(k)}(s-u)}\vc w_{k}\vc e_{j_{k}}^\top e^{(\ma T-\gamma \ma I)u} du\right\} (-\ma T)\vc 1\\\nonumber&=&C_k\left\{ \vc\alpha \vc u_{k}\vc e_{j_{k}}^\top [I-e^{(\ma T-\gamma \ma I)s}](\gamma \ma I-\ma T)^{-1}-\vc\alpha\mathcal{I}_{12}(s)-\vc v_k\mathcal{J}_{1k}(s)\right\}(-\ma T)\vc 1,
\end{eqnarray}where, by Lemma \ref{lemcarb}, $$\mathcal{I}_{12}(s)=(\vc f_{2,1}^\top\otimes I_m) e^{\ma B^{(k)}s}(\vc f_{2,2}\otimes I_m), \quad  \mathcal{J}_{1k}(s)=[\ma I_{mk},\ma 0_{mk\times m}]e^{\ma C^{(k)}s}[\ma 0_{m\times mk},\ma I_m]^\top$$ with 
$$\ma B^{(k)}=\left[\begin{array}{cc} \ma T &\vc u_{k}\vc e_{j_{k}}^\top\\ \ma 0 &\ma T-\gamma \ma I\end{array}\right],\quad \ma C^{(k)}=\left[\begin{array}{cc} \ma A^{(k)} &\vc w_{k}\vc e_{j_{k}}^\top\\ \ma 0 &\ma T-\gamma \ma I\end{array}\right].$$ 
\qed

\section{The uniform observation scheme} 
\begin{figure}[t] 
	\begin{center}
	\begin{picture}(300,150)(0,10)
	\includegraphics[angle=0,width=10cm]{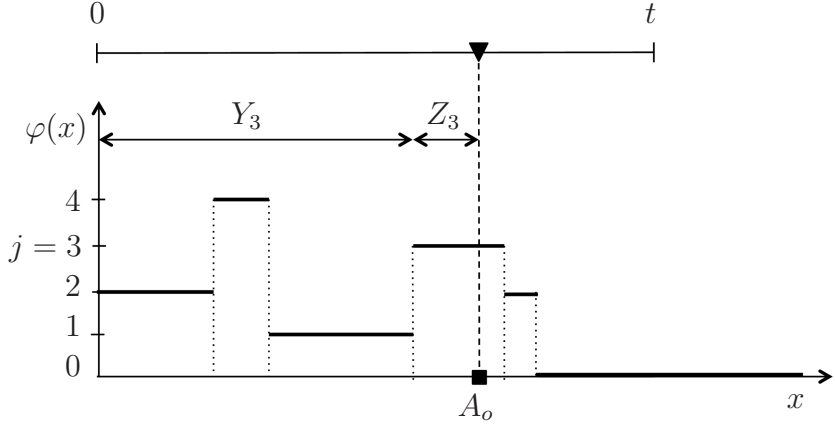} 
	\put(-283,138){0}\put(-75,138){$t$}
	\put(-292,5){0}
	\put(-292,20){$1$}\put(-292,35){2}\put(-313,50){$j=3$}\put(-292,67){4}\put(-230,100){$Y_3$}\put(-157,100){$Z_3$}\put(-145,-10){$A_o$}\put(-22,-7){$x$}\put(-307,95){$\varphi(x)$}
	\end{picture}\vspace*{1cm}
	\caption{\label{figpath2} A Possible trajectory of the phase process until absorption. The observation time is uniform on $[0,t]$ and is represented by a $ \blacktriangledown$ symbol. The observed phase is $j=3$.}
	\end{center}\end{figure}

In this section, we assume that an observer samples an individual in a population at a single random time $T_o$ in accordance with a uniform distribution on $[0,t]$, for some time $t>0$, where the clock is set at the birth of the individual. 
We then ask the same questions as in the Poisson observation scheme, but we expect different answers. The three random variables of interest $A_o, Y_j, $ and $Z_j$ are illustrated in Figure \ref{figpath2}. Their respective conditional distribution, under the uniform observation scheme, is provided in the next three propositions.

\begin{prop}\label{unifsingle} The conditional age distribution of the individual at a random observation time uniformly distributed on $[0,t]$, given $\varphi_o=j$, is given by 
\begin{equation}\label{qusol}\p[A_o\leq s\,|\,\varphi_o=j ]=1-\dfrac{ \vc\alpha [\exp(\ma Q s)-\exp(\ma Qt)]\,(-\ma Q)^{-1} \vc{e}_{j}}{\vc\alpha [\ma I-\exp(\ma Qt)](-\ma Q)^{-1} \vc{e}_j}\quad \mbox{for $s\leq t$,}\end{equation}
and $\p[A_o\leq s\,|\,\varphi_o=j ]=1$ for $s>t$.
\end{prop}

\noindent \textit{Proof.} By conditioning on the value of the observation time $T_o$, we have
\begin{equation}\label{po}\p[\varphi_o=j]=(1/t)\int_0^t \vc \alpha e^{\ma Qu}\vc e_j \, du=(1/t)\,\vc\alpha [\ma I-\exp(\ma Qt)](-\ma Q)^{-1} \vc{e}_j,\end{equation}and
\begin{eqnarray}\nonumber
\p[A_o\leq s, \varphi_o=j]&=&(1/t)\int_0^t \vc \alpha e^{\ma Qu}\vc e_j \mathds{1}_{\{u\leq s\}} du\\\nonumber&=&(1/t)\int_0^{\min(s,t)} \vc \alpha e^{\ma Qu}\vc e_j  du\\\label{p1}&=&(1/t)\, \vc\alpha [\ma I-\exp(\ma Q\min(s,t))]\,(-\ma Q)^{-1} \vc{e}_{j}.\end{eqnarray} The conditional distribution \eqref{qusol} then follows by dividing \eqref{p1} by \eqref{po}.
\qed
\medskip


\begin{prop} \label{pyju} The conditional distribution of $Y_j$, given $\varphi_o=j$, has a point mass at zero given by
\begin{eqnarray}\label{Yj0b} \p[Y_j=0\,|\,\varphi_o=j]&=&\dfrac{\alpha_j (1-e^{-\lambda_jt})}{\lambda_j\,\vc\alpha [\ma I-\exp(\ma Qt)](-\ma Q)^{-1} \vc{e}_j},\end{eqnarray} and for $0<y\leq t$,
\begin{eqnarray} \label{Yjb}\p[Y_j\leq y\,|\,\varphi_o=j]&=&1-\dfrac{\lambda_j\vc\alpha [e^{\ma Qy}-e^{\ma Q t}](-\ma Q)^{-1}\vc e_j +\vc\alpha e^{\ma Q y}\vc e_j(e^{-\lambda_j(t-y)}-1)}{\lambda_j\,\vc\alpha [\ma I-\exp(\ma Qt)](-\ma Q)^{-1} \vc{e}_j},\qquad\end{eqnarray}where $\lambda_j=-Q_{jj}$. Finally, $\p[Y_j\leq y\,|\,\varphi_o=j ]=1$ for $y>t$.
\end{prop}
\noindent \textit{Proof.} The proof follows exactly the same lines as in the Poisson observation case. For any $0<y\leq t$, 
\begin{eqnarray}\nonumber\lefteqn{\p[Y_j\in[y,y+dy], \,\varphi_o=j]}\\\nonumber&=&(1/t)\int_0^t\p[Y_j\in[y,y+dy], \,\varphi(u)=j\,|\,T_o\in[u,u+du]]\,du\\\nonumber&=& (1/t)\int_y^t\sum_{k\neq j}(\vc\alpha e^{\ma Q y})_k Q_{kj} dy\,e^{-\lambda_j (u-y)}\,du\\\nonumber&=& (1/t)\vc\alpha e^{\ma Q y} (\ma Q+\lambda_j\ma  I)\vc e_j\,dy\,\int_y^te^{-\lambda_j (u-y)}\,du\\\label{lasteq}&=& (1/t)\vc\alpha e^{\ma Q y} (\ma Q+\lambda_j \ma I)\vc e_j\,dy\,\dfrac{1- e^{-\lambda_j(t-y) }}{\lambda_j}.
\end{eqnarray}
Similarly, $$ \p[Y_j=0, \,\varphi_o=j]=(1/t)\int_y^t \alpha_j e^{-\lambda_j (u-y)}\,du=(1/t)\,\dfrac{\alpha_j (1- e^{-\lambda_j(t-y) })}{\lambda_j}, $$ which, together with \eqref{po}, leads to \eqref{Yj0b}. 
Finally, from \eqref{Yj0b}, \eqref{lasteq}, and \eqref{po} we obtain \eqref{Yjb}.
\qed
\medskip

\begin{prop}\label{pzju} The conditional distribution of $Z_j$, given $\varphi_o=j$, is given by
\begin{equation}\label{densZ1}\p[Z_j\leq z\,|\,\varphi_o=j]=1-\dfrac{\vc \alpha[\ma I-e^{\ma Q(t-z)}]\,(-\ma Q)^{-1}\vc e_j\,e^{-\lambda_j z} }{\vc\alpha  [\ma I-e^{\ma Qt}](-\ma Q)^{-1}\vc e_j}\quad\mbox{for $z\leq t$},\end{equation} and $\p[Z_j\leq z\,|\,\varphi_o=j]=1$ for $z>t$.
\end{prop}
\noindent \textit{Proof.} By the usual arguments,
\begin{eqnarray*}\lefteqn{\p[Z_j\in[z,z+dz],\,\varphi_o=j]}\\&=& (1/t)\left[\alpha_j e^{-\lambda_j z}+\int_z^t \sum_{k\neq j}(\vc\alpha e^{\ma Q(u-z)})_k Q_{kj} e^{-\lambda_j z}\,du\right] dz,
\end{eqnarray*} where the first term in the bracket accounts for the case where ``$T_o=z$'', that is, the individual is observed in her/his initial phase. We then have 
\begin{eqnarray}\nonumber\lefteqn{\p[Z_j\in[z,z+dz],\,\varphi_o	=j]}\\\nonumber&=&(1/t)\left[\alpha_j e^{-\lambda_j z}+\vc\alpha \int_z^{t} e^{\ma Q(u-z)} du\, (\ma Q+\lambda_j \ma I)\vc e_j e^{-\lambda_j z}\right]dz\\\label{Z2}&=&(1/t)\left[\alpha_j e^{-\lambda_j z}+\vc\alpha [\ma I-e^{\ma Q(t-z)}](- \ma Q)^{-1}\, (\ma Q+\lambda_j \ma I)\,\vc e_j e^{-\lambda_j z}\right]dz,\end{eqnarray}and the conditional density function of $Z_j$, given $\varphi_o=j$, is obtained by dividing \eqref{Z2} by \eqref{po} and rearranging the terms in the numerator:
\begin{equation*} f_{Z_j|\varphi_o=j}(z)=\dfrac{\vc\alpha \left[\lambda_j e^{-\lambda_j z}- e^{\ma Qt}e^{-(\ma Q+\lambda_j\ma I)z}(\ma Q+\lambda_j\ma I)\right] (-\ma Q)^{-1}\,\vc e_j}{\vc\alpha  [\ma I-e^{\ma Qt}](-\ma Q)^{-1}\vc e_j}.\end{equation*}Then, as $\p[Z_j\leq z\,|\,\varphi_o=j]=\int_0^z f_{Z_j|\varphi_o=j}(u)\,du$, we obtain \eqref{densZ1}. 
\qed


\section{Rare observation limit}\label{ro}
In practice, individuals of an animal population are usually observed very seldom. This is particularly true for endangered wild populations such as the Chatham Island black robins \textit{Petroica traversi}, which are observed once or twice per year (per individual) on average. We are therefore interested in the limit of the conditional age distribution as $\gamma\rightarrow0$ in the Poisson observation scheme, or as $t\rightarrow\infty$ in the uniform observation scheme.  First, observe that
$$\ma T\xrightarrow{\gamma\rightarrow0}\ma Q,\quad\mbox{and}\quad e^{\ma Qt}\xrightarrow{t\rightarrow\infty}\ma 0.$$ An interesting consequence of Propositions \ref{poissingle} and \ref{unifsingle} is that the rare observation limit of the age distribution is identical for the Poisson and uniform observation schemes. In addition, the limiting age distribution corresponds to the age distribution $F_j(s)$ that we derived in Lemma \ref{pbirth} under the assumption that the birth process is Poisson.



\begin{cor} \label{rareobs}The rare observation limit of the conditional age distribution, given $\varphi_o=j$, is given by\begin{eqnarray}\label{q5sol}\p[A_o\leq s\,|\,\varphi_o=j]&=& 1-\dfrac{ \vc\alpha \exp(\ma Q s)\,(-\ma Q)^{-1} \vc e_j}{\vc\alpha (-\ma Q)^{-1} \vc e_j}=F_j(s).
\end{eqnarray}
\qed\end{cor}

Actually, not only is the rare observation limit of the age distribution identical for the Poisson and the uniform observation schemes, but this holds for the conditional distribution of $Y_j$ and $Z_j$ too, as shown in the next two corollaries. These results are direct consequences of Propositions \ref{pyjp} and \ref{pyju}, and Propositions \ref{pzjp} and \ref{pzju}, respectively. 

\begin{cor}The rare observation limit of the conditional distribution of $Y_j$, given $\varphi_o=j$,
has a point mass at zero and is given by
\begin{eqnarray}\label{Yj0} \p[Y_j=0\,|\,\varphi_o=j]&=&\dfrac{\alpha_j}{\lambda_j\,\vc\alpha (-\ma Q)^{-1}\vc e_j}\end{eqnarray} and for $y>0$,\begin{eqnarray} \label{Yj}\p[Y_j\leq y\,|\,\varphi_o=j]&=&1-\dfrac{\vc\alpha e^{{\bf Q}y}(-{\bf Q})^{-1}({\bf Q}+\lambda_j \,{\bf I})\vc e_j}{\lambda_j\,\vc\alpha (-{\bf Q})^{-1}\vc e_j},\end{eqnarray}where $\lambda_j=-Q_{jj}$.\qed
\end{cor}

\begin{cor} The rare observation limit of the conditional distribution of $Z_j$, given $\varphi_o=j$, is exponential with parameter $\lambda_j$. \qed
\end{cor}

\section{Numerical illustrations}

We illustrate the results of the previous sections on a toy example first, and then on the real-world example of the Chatham Island black robin \textit{Petroica traversi} population.
\subsection{Toy example with five phases}

We consider a PH$(\vc\alpha,Q)$ lifetime distribution with $m=5$ transient phases and transition rate matrix $$Q=\left[\begin{array}{ccccc}-3& 2& &&\\
    & -5& 3& &\\
    1&& -4 &2& \\
     &1& & -6 &3\\
     &&1 & & -2\end{array}\right],$$ where a blank space represents a zero entry. We shall consider two initial distribution vectors: $$\vc\alpha^{(1)}=[1,0,0,0,0],\quad\mbox{and}\quad\vc\alpha^{(2)}=[1/5,1/5,1/5,1/5,1/5].$$ In the first case, the process starts in phase 1 almost surely, while in the second case, the initial phase is chosen uniformly on the transient phase space. As we show below, the initial distribution vector can affect the shape of the various conditional distributions. We choose to represent densities (rather than distribution functions) as they better capture the features of the distributions.
     
\begin{figure}[t] 
	\begin{center}
	\includegraphics[angle=0,width=10cm]{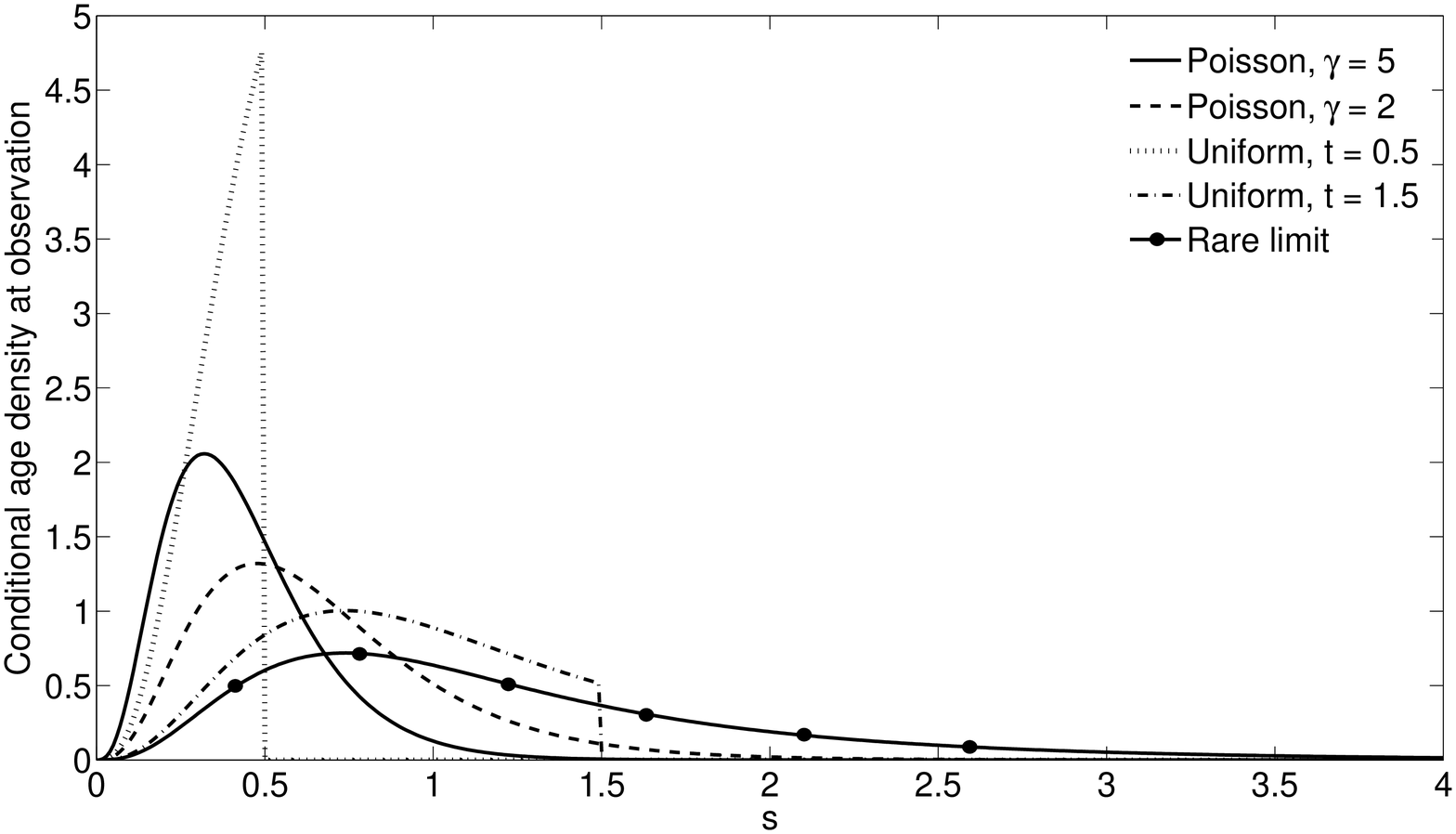} 
	\includegraphics[angle=0,width=10cm]{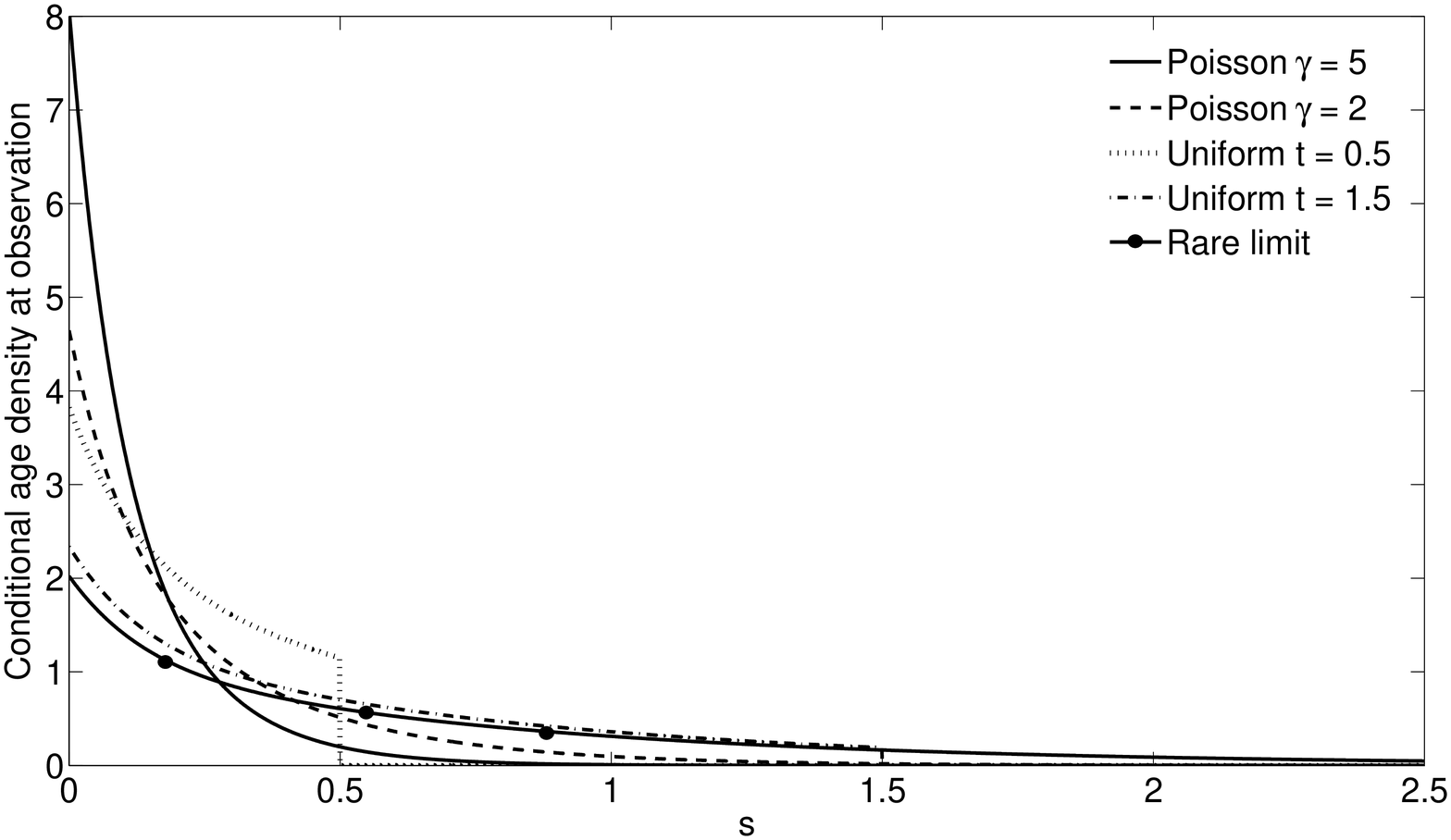} 
	\caption{\label{fig1}Conditional age density at the observation time, given that phase $j=4$ is observed, for the two observation schemes with different parameter values, and the rare observation limit, as $\vc\alpha=[1,0,0,0,0]$ (top) and $\vc\alpha=[1/5,1/5,1/5,1/5,1/5]$ (bottom). }
	\end{center}\end{figure}

  We first assume a single observation, and we condition on the observed phase being $j=4$. Figure \ref{fig1} shows the conditional age densities obtained under the different observation schemes, for different values of the parameters $\gamma$ and $t$, as well as the rare observation limit. Observe the discontinuity of the density at $s=t$ in the uniform case (while the distribution function given in \eqref{qusol} is continuous at $s=t$). This comes from the fact that an individual observed in the time window $[0,t]$ cannot be older than $t$, and suggests that the Poisson observation scheme is more natural than the uniform observation scheme on a finite time interval.
 We see that for $\vc\alpha^{(1)}$, the mode of the distribution is positive and tends to increase as the observation becomes rare, while for $\vc\alpha^{(2)}$, the mode is clearly at age zero. 

\begin{figure}[t] 
	\begin{center}
		\includegraphics[angle=0,width=10cm]{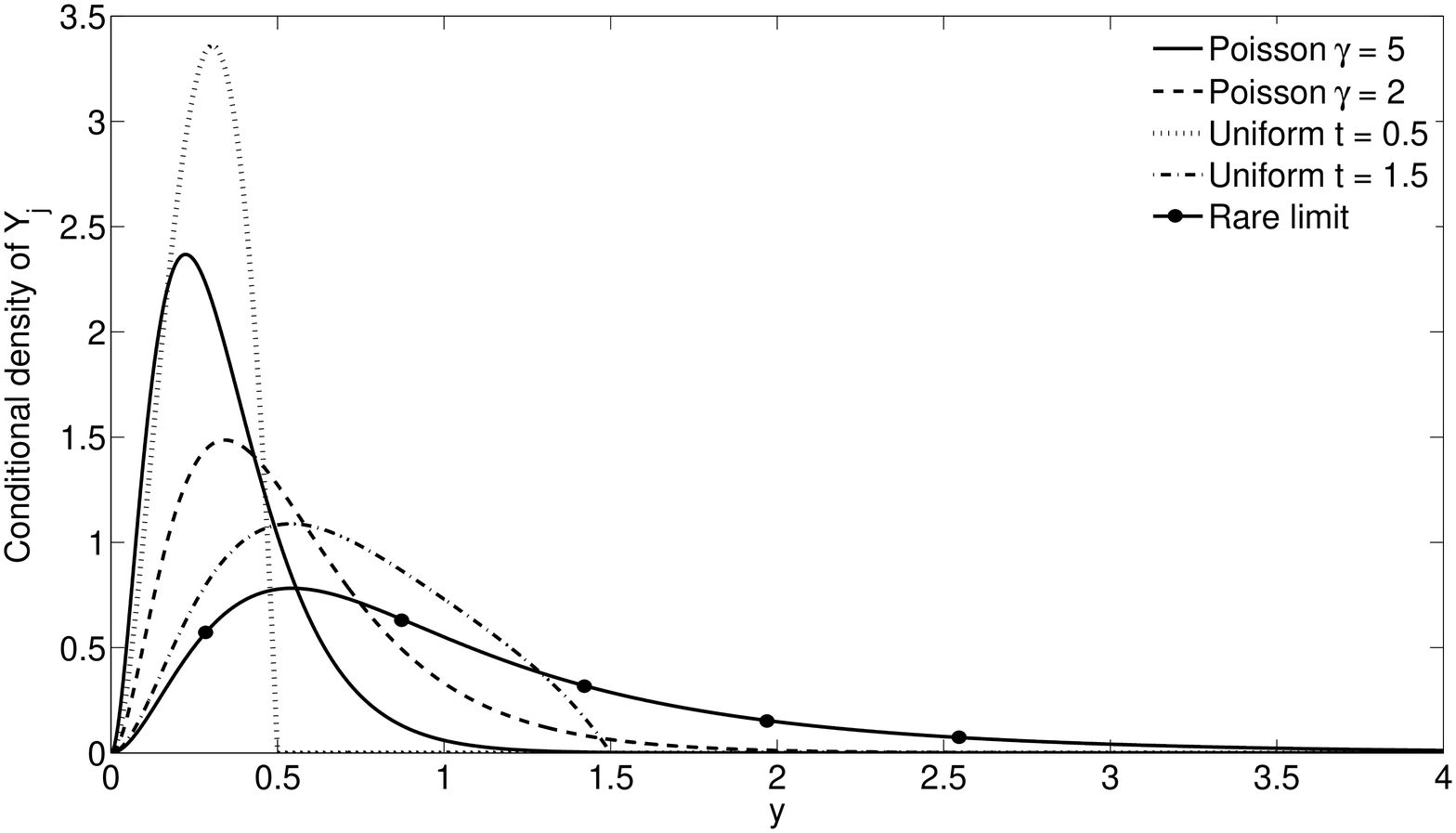} 
	\includegraphics[angle=0,width=10cm]{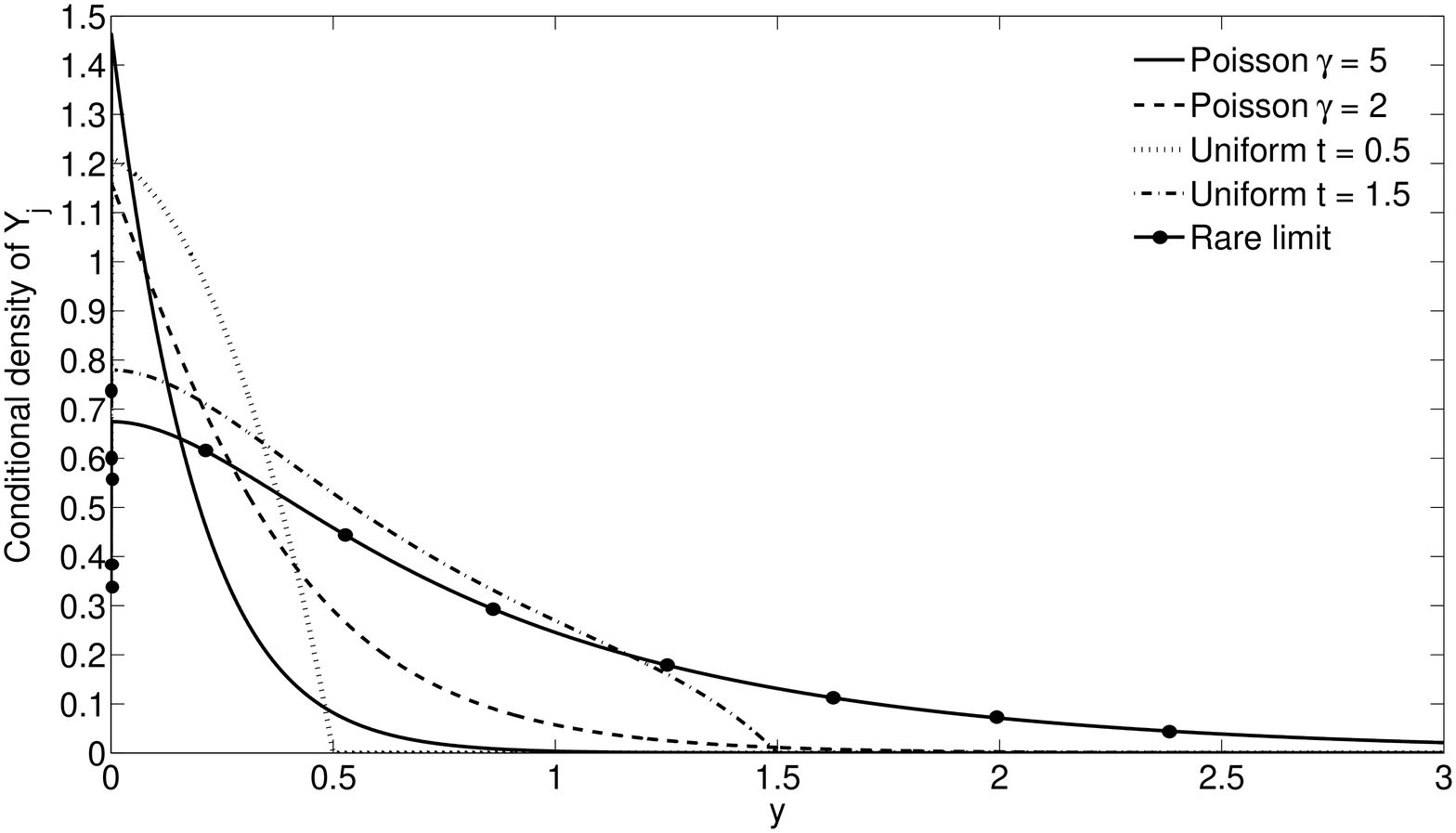} 
	\caption{\label{fig1b}Conditional density of $Y_j$, given that phase $j=4$ is observed, for the two observation schemes with different parameter values, and the rare observation limit, as $\vc\alpha=[1,0,0,0,0]$ (top) and $\vc\alpha=[1/5,1/5,1/5,1/5,1/5]$ (bottom). The point mass at zero is clear when $\vc\alpha=[1/5,1/5,1/5,1/5,1/5]$.}
	\end{center}\end{figure}

\begin{figure}[t] 
	\begin{center}
		\includegraphics[angle=0,width=10cm]{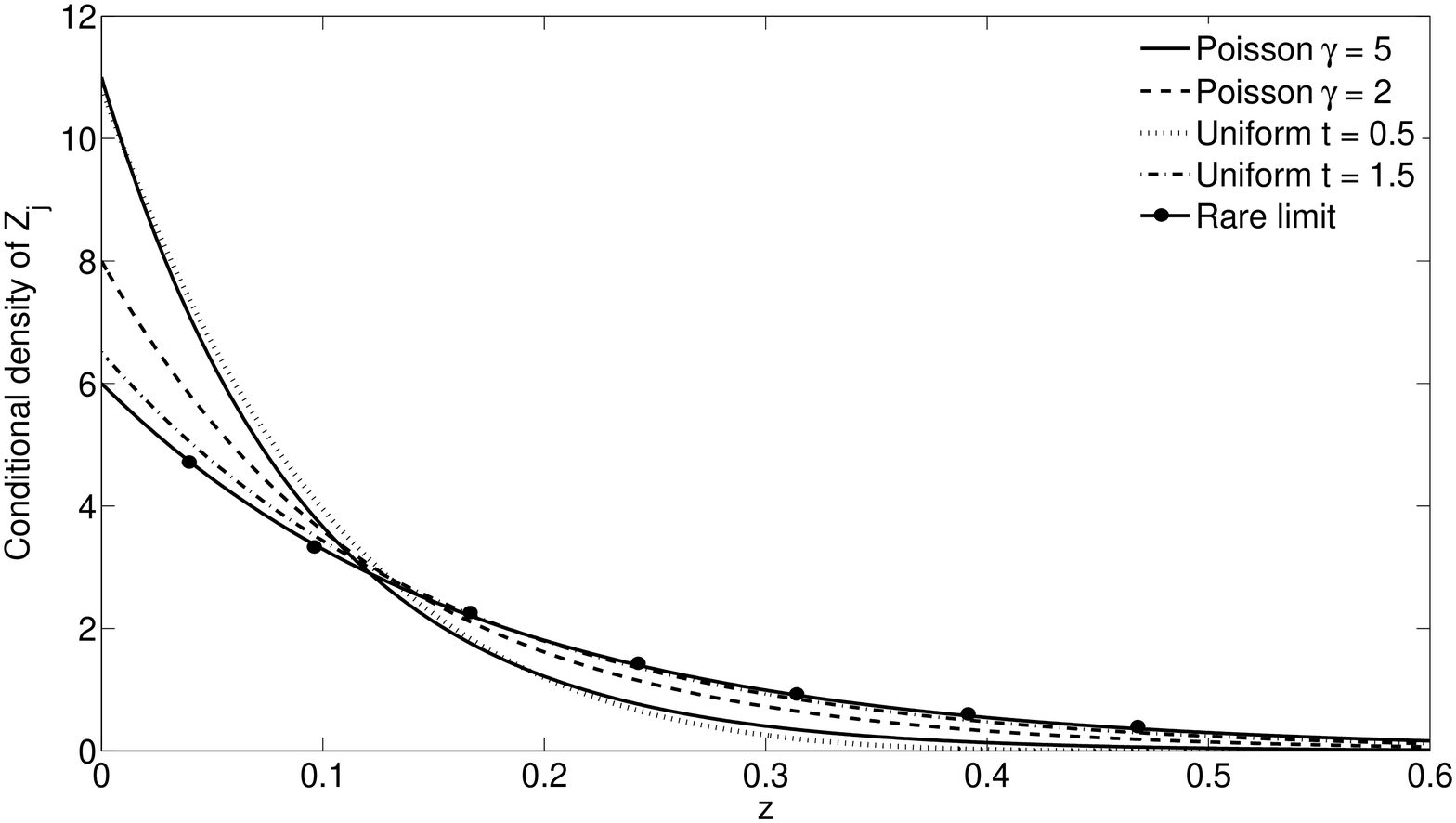} 
	\includegraphics[angle=0,width=10cm]{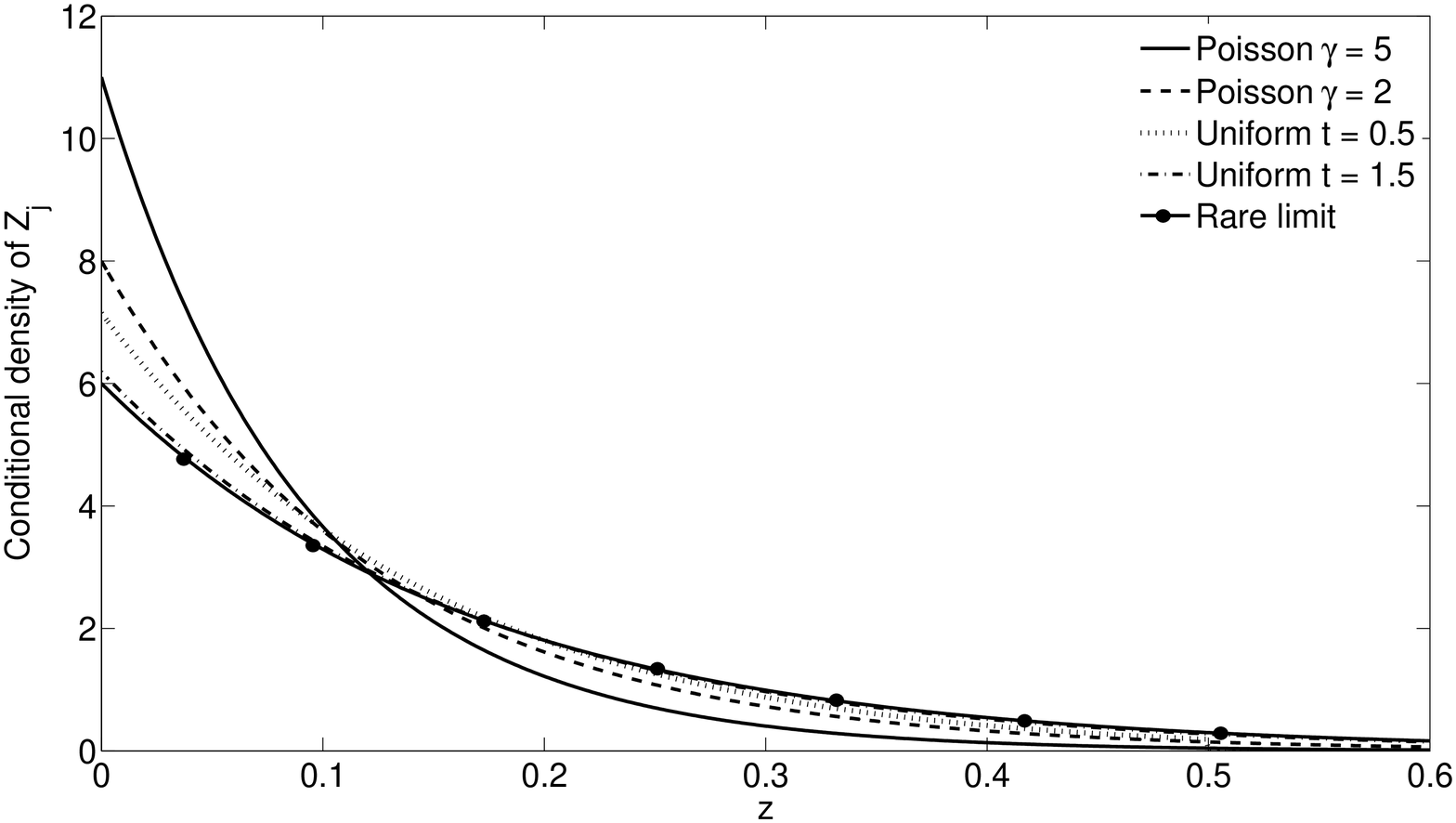} 
	\caption{\label{fig1c}Conditional density of $Z_j$, given that phase $j=4$ is observed, for the two observation schemes with different parameter values, and the rare observation limit, as $\vc\alpha=[1,0,0,0,0]$ (top) and $\vc\alpha=[1/5,1/5,1/5,1/5,1/5]$ (bottom).}
	\end{center}\end{figure}

Figures \ref{fig1b} and \ref{fig1c} illustrate the conditional densities of $Y_4$ and $Z_4$ respectively. The initial distribution affects the shape of the density of $Y_4$, which has a point mass at zero for $\vc\alpha^{(2)}$, but has negligible effect on the shape of the density of $Z_4$. Also note that in the uniform observation case, the density of $Y_j$ is continuous at $s=t$, but similar to the age density in Figure \ref{fig1}, the density of $Z_j$ is discontinuous at $s=t$.

\begin{figure}[t] 
	\begin{center}
		\includegraphics[angle=0,width=10cm]{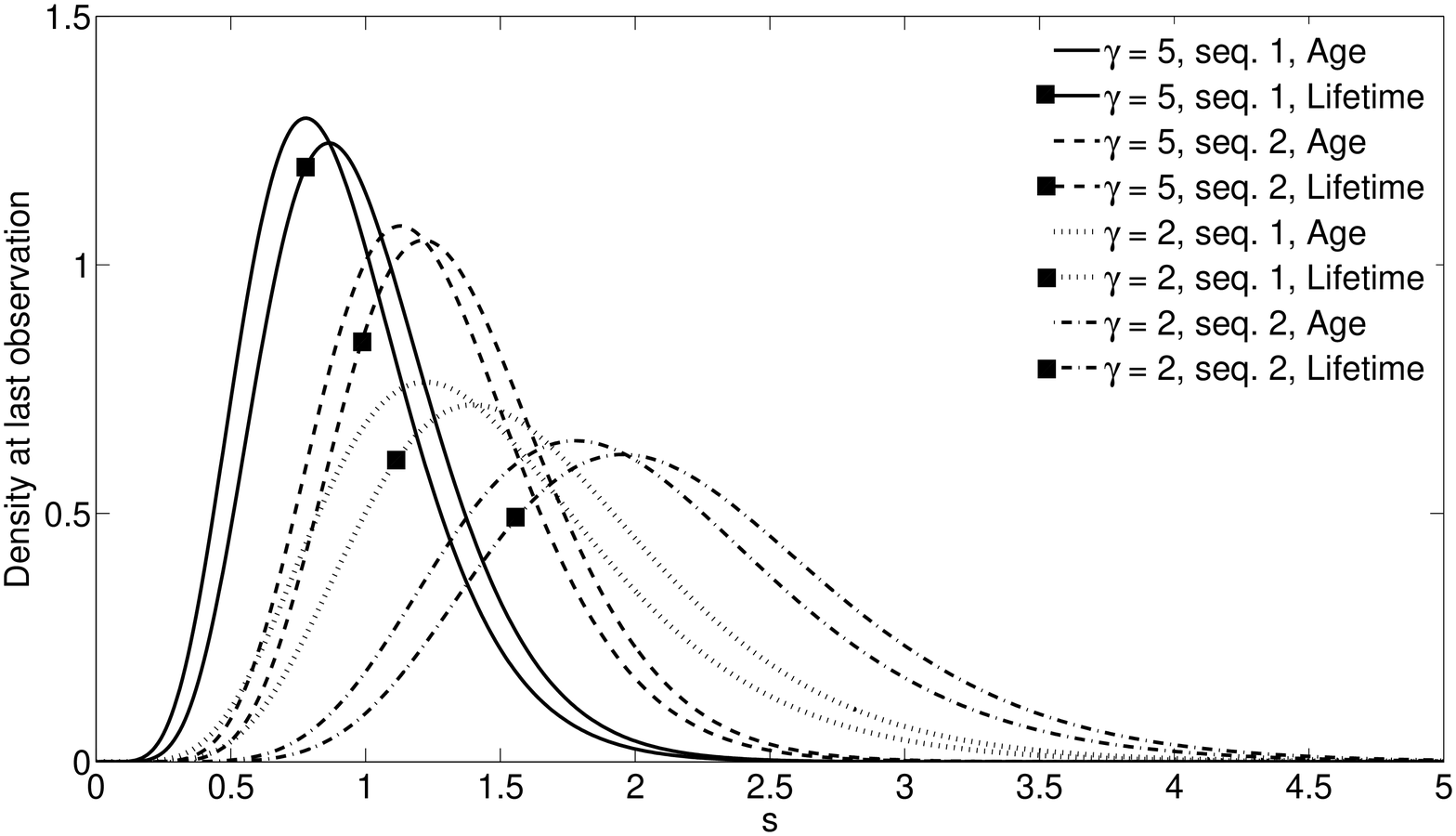} 
	\caption{\label{fig2}Conditional age/lifetime density at the time of the last Poisson observation when the observed sequence is $1,2,3,3,4$ (Sequence 1) and $2,3,1,1,4$ (Sequence 2), with $\vc\alpha=[1,0,0,0,0]$.}
	\end{center}\end{figure}

Finally, we consider the multiple Poisson observation scheme with $k=5$ observations, for different values of the parameter $\gamma$, and different sequences of observed phases: $1,2,3,3,4$ (Sequence 1) and $2,3,1,1,4$ (Sequence 2, which is less likely than Sequence 1). In Figure \ref{fig2} we compare the age distribution at the $5$th observation, and the lifetime distribution given that the individual is dead at the $6$th observation, for $\vc\alpha^{(1)}$ (there is not much difference for $\vc\alpha^{(2)}$). 
The graph illustrates how a change in the sequence of observed phases affects the related conditional distributions. We see that the tail of the distributions corresponding to Sequence 2 is fatter than for Sequence 1, that is, an individual is more likely to be older at the time of the last observation when Sequence 2 is observed.

	\subsection{How old are the Chatham Island black robins \textit{Petroica traversi}?}
	
	In this last section, we come back to our original objective, and illustrate the usefulness of our results to compute the age pyramid for the black robin population during the intensive management period between 1980 and 1989. 

A first step of the analysis consists in modelling the bird population using a branching process called Markovian binary tree (MBT), which is done in detail in a parallel study\footnote{S. Hautphenne, M. Massaro, E. S. Kennedy, and R. Sainudiin. Modelling of the Chatham Island black robin \textit{Petroica traversi} populations using branching processes: Informed management strategies for reintroduction of endangered species. In preparation}. Age-specific mortality and fertility rates of the black robins can be estimated from the unique dataset collected between 1980 and 1989 \cite{butlerblack}. These age-specific rates are used to estimate the parameters of an MBT that optimally fits the data, and this model is then used to study demographic properties of the population during the intensive management period. 

In the present section we shall focus on the bird lifetime distribution rather than on their reproduction process. The estimated female age-specific mortality rates for the period 1980-1989 are shown in Table~\ref{t1}. The age class $[0,1)$ corresponds to birds who fledged. Note that the lack of data, in particular for the ages above 5, lead to inacurrate estimations for these age-classes. Indeed, five females reached age 5 between 1980 and 1989, but only one reached age 8, and only one reached age 12 during that period. A Bayesian approach was used to bias low relative frequencies upward.  

\begin{figure}[t] 
	\begin{center}
	\includegraphics[angle=0,width=10cm]{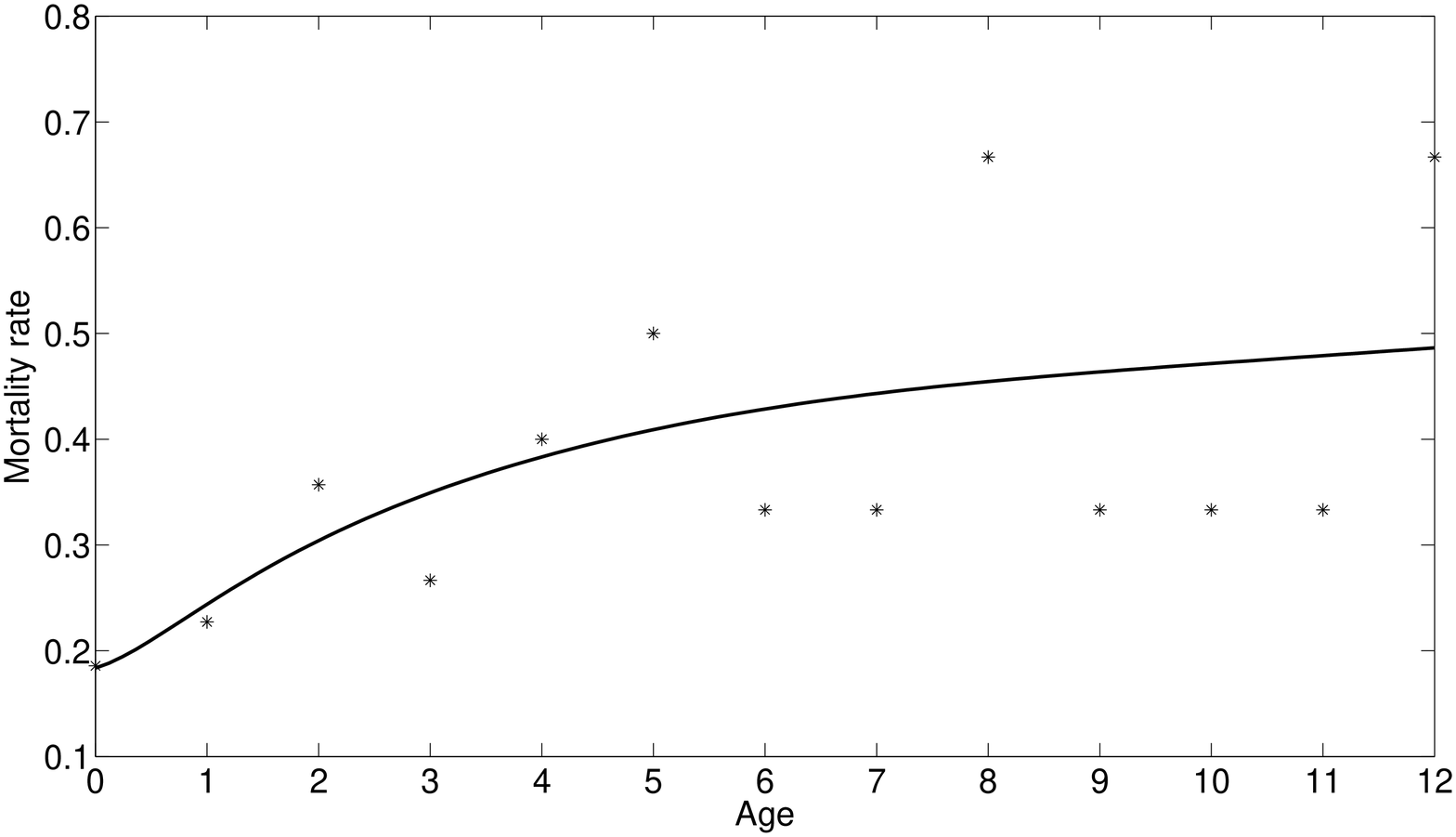} 
	\caption{\label{figdon} Estimated age-specific mortality rates (stars) and the age-specific mortality curve computed with the optimal Markovian model fitting the data (plain line).}
	\end{center}\end{figure}

We assume that the lifetime $L$ of a female black robin has a PH$(\vc\alpha,Q)$ distribution with $m=13$ (transient) phases with the specific ageing structure
\begin{equation}\label{structure} \vc\alpha=\kbordermatrix{\mbox{}&1&2&\ldots&m\\ &1&0&\ldots&0}, \quad
Q=\kbordermatrix{\mbox{}&1&2&3&\ldots&m\\
1&-\lambda_1&\lambda_1s_1&& &&\\
2&&-\lambda_2&\lambda_2s_2&&\\
\vdots& && &\ddots&\\
m&&& &&-\lambda_m\\
};\end{equation} that is, an individual starts its life in phase 1 and moves through successive phases until it dies; in this case, a transition from a transient phase $j$ can only be to the next phase $j+1$ (with probability $s_j$), or to the absorbing phase 0 (with probabilty $1-s_j$). This particular PH distribution has a minimal number of parameters ($2m-1$), which makes it more convenient for parameter estimation. The rates $\lambda_i$ and probabilities $s_i$ were estimated following the approach in \cite{lin2007markov}, by minimizing the sum of weighted squared errors
$$F=\sum_{x=0}^{12}(\hat{d}_x-\bar{d}(x))^2 \hat{S}_x,$$ where $\hat{d}_x$ is the observed (estimated) mortality rate in age class $[x,x+1)$ as shown in Table \ref{t1}, $\hat{S}_x=(1-\hat{d}_0)(1-\hat{d}_1)\cdots (1-\hat{d}_{x-1})$ is the observed probability of survival until age class $[x,x+1)$, and $\bar{d}(x)$ is the corresponding model value for $\hat{d}_x$, which is computed using \eqref{disPH}:
\begin{eqnarray*} \bar{d}(x) &=&\p[x< L \leq x+1 | L> x]\\ &=& \dfrac{\p[L>x]-\p[L> x+1]}{\p[L>x]}\\&=&\dfrac{\vc\alpha e^{Q x}(I-e^Q)\vc 1}{\vc\alpha e^{Q x}\vc 1}.
\end{eqnarray*}
The resulting optimal rates $\lambda_i$ and probabilities $s_i$ are provided in Table~\ref{t2}. Figure \ref{figdon} shows the estimated age-specific mortality rates $\hat{d}_x$ (stars) together with the mortality function $\bar{d}(x)$ corresponding to the optimal model (plain line). We see that the model mortality curve smoothes the inaccurate point estimates in a satisfactory way.


\begin{table}\centering\begin{tabular}{ll} Age class $[x,x+1)$ & Mortality rate $\hat{d}_x$\\\hline $[0,1)$&0.19\\$[1,2)$& 0.23 \\$[2,3)$&0.36\\$[3,4)$& 0.27\\$[4,5)$& 0.40 \\$[5,6)$& 0.50\\$[6,7)$& 0.33\\ $[7,8)$& 0.33\\$[8,9)$& 0.67\\$[9,10)$& 0.33\\$[10,11)$& 0.33\\$[11,12)$& 0.33\\$[12,\infty)$& 0.67\end{tabular}\caption{\label{t1}Female age-specific mortality rates estimated from the raw data on the Chatham Island black robin population during the managed phase 1980-1989.}\end{table} 

\begin{table}\centering\begin{tabular}{lll}  Phase $i$ & $\lambda_i$ & $s_i$\\\hline 1 & 1.54 & 0.93\\2 & 1.59 & 1 \\3 & 1.26& 0.27\\4&1.53& 1\\5& 1.67 & 1 \\6& 2.21 & 0.98\\7& 1.86 & 1\\ 8 &1.36 &1\\9 &1.28 &0.06\\10 & 0.68&0.98\\11 &0.76 &0.91\\12 & 1.86&0.55\\13&1.15&$-$\end{tabular}\caption{\label{t2}The parameters of the PH distribution modelling the  Chatham Island black robin lifetime.}\end{table} 

\begin{figure}[t] 
	\begin{center}
	\includegraphics[angle=0,width=10cm]{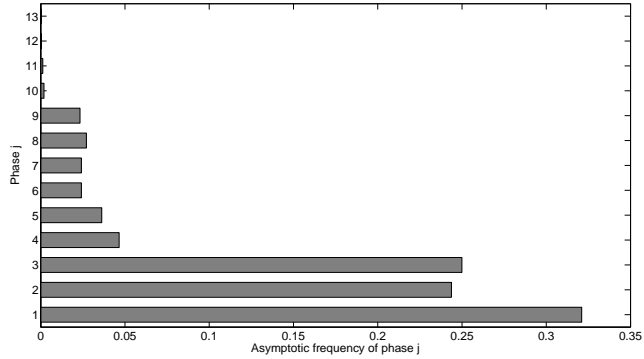} 
	\caption{\label{fig3}The asymptotic phase frequency in the Chatham Island black robin population if the population evolves with the  demographic rates of 1980-1989.}
	\end{center}\end{figure}

Among other useful properties, the MBT model allows us to compute the asymptotic phase frequency in the population, that is, the proportion of birds in each of the 13 phases if we let the population evolve for a long period of time with the same demographic rates. We denote by $f^p_j$ the asymptotic frequency of phase $j$, and we show the value of $f^p_j$ for $1\leq j\leq 13$ in Figure~\ref{fig3}.
	Since the phases do not have any physical interpretation, the asymptotic phase frequency does not have much biological interest in its own. However, it can be used in combinaison with the results developed in this paper to compute the asymptotic age-frequency  (also called the age-pyramid), which cannot be obtained directly from the MBT model.
	
	We consider the following eleven age-classes: $[0,1)$, $[1,2)$,\ldots,$[9,10)$, $[10,\infty)$, and we denote the asymptotic frequency of the age-class starting at age $x$ as $f^a_x$, $0\leq x\leq 10$. We can then compute $f^a_x$ as
	$$f^a_x=\sum_{1\leq j\leq 13} f^p_j \,\p[\mbox{age}\in [x,x+1)\,|\, \mbox{phase} = j].$$ We approximate the probability $\p[\mbox{age}\in [x,x+1)\,|\, \mbox{phase} = j]$ using the rare limit conditional age distribution provided in Corollary \ref{rareobs}. 
	
\begin{figure}[t] 
	\begin{center}
	\includegraphics[angle=0,width=10cm]{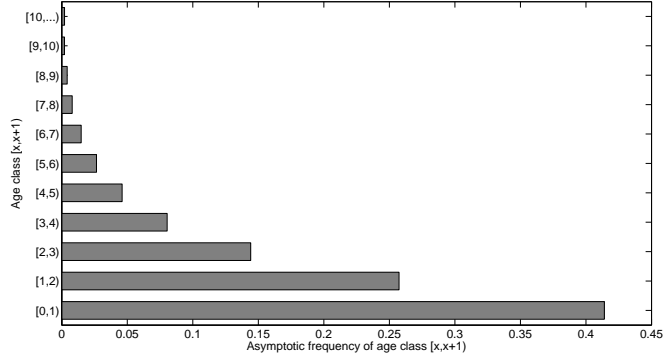} 
	\caption{\label{fig4}The asymptotic age frequency in the Chatham Island black robin population if the population evolves with the  demographic rates of 1980-1989.}
	\end{center}\end{figure}

	The resulting age pyramid for the black robins is depicted in Figure~\ref{fig4}. From the shape of the pyramid, we see that the population is rapidly expanding, which was indeed the case during the period of intensive conservation management in 1980-1989.
	However, during this period, reproductive outputs were artificially increased through human intervention. By cross-fostering black robin offspring to the closely related Chatham Island tomtit \textit{Petroica macrocephala chathamensis}, female black robins were induced to lay additional clutches of eggs \cite{massaro2013human}. This artificially increased reproductive success in combination with our assumption that the fertility and mortality rates are fixed over a long time period, results in an age pyramid whose shape may not be representative of the current population.
	
	 Given that the current black robin population is restricted to only two small islands and includes fewer than 250 individuals, the species remains endangered \cite{massaro2013nest}. Hence, knowing the age frequency of this population, and the associated fertility and mortality rates, is highly relevant to the future conservation management of this species; this is investigated in more detail
	  in the parallel study\footnote{\label{ff}S. Hautphenne, M. Massaro, E. S. Kennedy and R. Sainudiin. Modelling of the Chatham Island black robin \textit{Petroica traversi} populations using branching processes: Informed management strategies for reintroduction of endangered species. In preparation}.
	
	\begin{rem}
	In the particular case of an ageing process with structure \eqref{structure}, and conditionally on phase $j$ being observed, any trajectory of the phase process before observation is restricted to the phases $1,2,\ldots,j$. Therefore, for $j<m$, it is sufficient to consider the process restricted to the smaller phase space $\{1,2,\ldots, j\}$, with initial distribution vector and generator
$$  \vc\alpha^{(j)}=\kbordermatrix{\mbox{}&1&2&\ldots&j\\ &1&0&\ldots&0}, \;
Q^{(j)}=\kbordermatrix{\mbox{}&1&2&3&\ldots&j\\
1&-\lambda_1&\lambda_1s_1&& &&\\
2&&-\lambda_2&\lambda_2s_2&&\\
\vdots& && &\ddots&\\
j&&& &&-\lambda_j\\
}.$$ 
We further define the matrix $T^{(j)}=Q^{(j)}-\gamma I$ where the identity matrix is $j\times j$, and $\vc e_i^{(j)}$ as the truncation of $\vc e_i$ after its $j$th entry, for $i\leq j$. The matrix ingredients $\vc\alpha, Q, T$, and $ \vc e_j$ used in the lemmas, propositions, and corollaries in the previous sections can then be replaced by their smaller counterpart $\vc \alpha^{(j)}$, $Q^{(j)}$, $T^{(j)}$, and $\vc e_{j}^{(j)}$.

Finally, observe that with the particular ageing structure \eqref{structure}, in the Poisson observation scheme, the random variable $Y_j$ has the same distribution as the time until the process with generator $T^{(j)}$ reaches phase $j$ for the first time, $B(j)$, conditionally on $B(j)<\infty$.
\end{rem}

\section*{acknowledgements}
The authors are supported by the Australian Research Council Laureate Fellowship FL130100039. The first author has also conducted part of the work under the Discovery Early Career Researcher Award DE150101044. Finally, we thank D. Merton, E. Kennedy, R. Morris, A. Munn, G. Murman, R. Nilsson, R. Thorpe and many other Wildlife Service and Department of Conservation staff that helped over the past 30 years to bring the black robin back from the brink of extinction.



\end{document}